\documentclass[10pt]{article}
\usepackage{mathtools}
\usepackage{graphicx}
\usepackage{amsmath,amsfonts,latexsym,amscd,amssymb,theorem}
\usepackage{graphicx}
\usepackage[ansinew]{inputenc}
\usepackage{epsfig}
\usepackage{amsmath}
\usepackage{amssymb}
\usepackage{afterpage}
\usepackage{selinput}
\usepackage{fancyhdr}
\newcommand{\ba}{\begin{eqnarray}}
\newcommand{\ea}{\end{eqnarray}}

\newtheorem{thm}{Theorem}[section]
\newtheorem{problem}{Problem}
\newtheorem{conjecture}{Conjecture}

\newtheorem{theorem}[thm]{Theorem}

\newtheorem{lemma}[thm]{Lemma}

\makeatletter
\newcommand*{\rom}[1]{\expandafter\@slowromancap\romannumeral #1@}
\makeatother

\begin{document}
\title{\textbf{Forbidding Couples of Tournaments and the Erd\"{o}s-Hajnal Conjecture}}
\maketitle


\begin{center}
\author{Soukaina ZAYAT \footnote{Department of Mathematics Faculty of Sciences I, Lebanese University, KALMA Laboratory, Beirut - Lebanon.\\ (soukaina.zayat.96@outlook.com)}, Salman GHAZAL \footnote{Department of Mathematics Faculty of Sciences I, Lebanese University, Beirut - Lebanon.\\ (salman.ghazal@ul.edu.lb)}\footnote{Department of Mathematics and Physics, School of Arts and Sciences, Beirut International University, Beirut - Lebanon. \vspace{2mm} (salman.ghazal@liu.edu.lb)}}
\end{center}

\begin{abstract}
A celebrated unresolved conjecture of Erd\"{o}s and Hajnal states that for every undirected graph $H$ there exists $ \epsilon(H) > 0 $ such that every undirected graph on $ n $ vertices that does not contain $H$ as an induced subgraph contains a clique or a stable set of size at least $ n^{\epsilon(H)} $. This conjecture has a directed equivalent version stating that for every tournament $H$ there exists $ \epsilon(H) > 0 $ such that every $H-$free $n-$vertex tournament $T$ contains a transitive subtournament of size at least $ n^{\epsilon(H)} $. Recently the conjecture was proved for all six-vertex tournaments, except $K_{6}$. In this paper we construct two infinite families of tournaments for which the conjecture is still open for infinitely many tournaments in these two families  $-$ the family of so-called super nebulas and the family of so-called super triangular galaxies. We prove that for every super nebula $H_{1}$ and every $\Delta$galaxy $H_{2}$ there exist $\epsilon(H_{1},H_{2})$ such that every $\lbrace H_{1},H_{2}\rbrace$$-$free tournament $T$ contains a transitive subtournament of size at least $\mid$$T$$\mid^{\epsilon(H_{1},H_{2})}$.  We also prove that for every central triangular galaxy $H$ there exist $\epsilon(K_{6},H)$ such that every $\lbrace K_{6},H\rbrace$$-$free tournament $T$ contains a transitive subtournament of size at least $\mid$$T$$\mid^{\epsilon(K_{6},H)}$. And we give an extension of our results.  \end{abstract}

\section{Introduction}
Let $ G $ be an undirected graph. We denote by $ V(G) $ the set of its vertices and by $ E(G) $ the set of its edges. We call $ \mid$$G$$\mid$ $=$ $ \mid$$V(G)$$\mid$ the \textit{size} of $G$. A \textit{clique} in $G$ is a set of pairwise adjacent vertices and a \textit{stable set} in $G$ is a set of pairwise nonadjacent vertices. A \textit{digraph} is a pair $D=(V,E)$ of sets such that $E\subset V \times V$, and such that for every $(x,y)\in E$ we must have $(y,x)\notin E$, in particular if $(x,y)\in E$ then $x \neq y$. $E$ is the arc set and $V$ is the vertex set and they are denoted by $E(D)$ and $V(D)$ respectively. We say that $D'$ is a \textit{subdigraph} of a digraph $D$ if $V(D') \subseteq V(D)$ and $E(D') \subseteq E(D)$. Let $X \subseteq V(D)$, the \textit{subdigraph of} $D$ \textit{induced by} $X$ is denoted by $D$$\mid$$X$, that is the digraph with vertex set $X$, such that for $x,y \in X$, $(x,y) \in E(D$$\mid$$ X)$ if and only if $(x,y) \in E(D)$. Denote by $D\backslash X$ the subdigraph of $D$ induced by $V(D)\backslash X$. We say that $D$ \textit{contains} $D'$ if $D'$ is isomorphic to a subdigraph of $D$. A \textit{tournament} is a directed graph (digraph) such that for every pair $u$ and $v$ of vertices, exactly one of the arcs $(u,v)$ or $(v,u)$ exists. A tournament is \textit{transitive} if it contains no directed cycle. Let $T$ be a tournament. We write $\mid$$T$$\mid$ for $\mid$$V(T)$$\mid$ and we say that $\mid$$T$$\mid$ is the \textit{size} of $T$. If $(u,v)\in E(T)$ then we say that $u$ is \textit{adjacent to} $v$ (alternatively: $v$ is an \textit{outneighbor} of $u$) and we write $u\rightarrow v$, also we say that $v$ is \textit{adjacent from} $u$ (alternatively: $u$ is an \textit{inneighbor} of $v$) and we write $v\leftarrow u$. For two sets of vertices $V_{1},V_{2}$ of $T$ we say that $V_{1}$ is \textit{complete to} (resp. \textit{from}) $V_{2}$ if every vertex of $V_{1}$ is adjacent to (resp. from) every vertex of $V_{2}$, and we write $V_{1} \rightarrow V_{2}$ (resp. $V_{1} \leftarrow V_{2}$). We say that a vertex $v$ is complete to (resp. from) a set $V$ if $\lbrace v \rbrace$ is complete to (resp. from) $V$ and we write $v \rightarrow V$ (resp. $v \leftarrow V$). Given a tournament $H$, we say that $T$ \textit{contains} $H$ if $H$ is isomorphic to $T$$\mid$$X$ for some $X \subseteq V(T)$. If $T$ does not contain $H$, we say that $T$ is $H$$-$$free$. For a class $\mathcal{F}$ of tournaments and a tournament $T$ we say that $T$ is $\mathcal{F}$$-$$free$ if $T$ is $F$$-$$free$ for every $F\in \mathcal{F}$. 
\vspace{1.5mm}\\
Erd\"{o}s and Hajnal proposed the following conjecture (EHC)\cite{jhp}:
\begin{conjecture} For any undirected graph $H$ there exists $ \epsilon(H) > 0 $ such that every $n$$-$vertex undirected graph that does not contain $H$ as an induced subgraph contains a clique or a stable set of size at least $ n^{\epsilon(H)}. $
\end{conjecture}
In 2001 Alon et al. proved \cite{fdo} that Conjecture $1$ has an equivalent directed version, as follows:
\begin{conjecture} \label{a} For any tournament $H$ there exists $ \epsilon(H) > 0 $ such that every $ H- $free tournament with $n$ vertices contains a transitive subtournament of size at least $ n^{\epsilon(H)}. $
\end{conjecture}
A class of tournaments $\mathcal{F}$ \textit{satisfy the Erd\"{o}s-Hajnal Conjecture (EHC)} (equivalently: $\mathcal{F}$ has the \textit{Erd\"{o}s-Hajnal property}) if there exists $ \epsilon(\mathcal{F}) > 0 $ such that every $\mathcal{F}$$- $free tournament $T$ with $n$ vertices contains a transitive subtournament of size at least $ n^{\epsilon(\mathcal{F})}. $ If $\lbrace H\rbrace$ satisfy $EHC$ we simply say that $H$ \textit{satisfies EHC}.\vspace{2mm}

Instead of forbidding just one tournament, one can state the analogous conjecture where we forbid two tournaments. The only results in this setting are in \cite{kgg,ssss}.\\
Conjecture \ref{a} has not yet been proved for \textit{super nebulas}, $K_{6}$, and \textit{super triangular galaxies}. That motivates the work of this paper. In this paper we prove that $\lbrace H_{1},H_{2}\rbrace$$-$free tournaments $T$ contain transitive subtournaments of size at least $\mid$$T$$\mid$$^{\epsilon (H_{1},H_{2})}$ for some $\epsilon (H_{1},H_{2})>0$ and infinite number of couples of tournaments $\lbrace H_{1},H_{2}\rbrace$.
Before stating formally our results, we need to introduce some definitions and notations.\vspace{2mm}

 Let $ \theta = (v_{1},...,v_{n}) $ be an ordering of the vertex set $V(D)$ of an $ n- $vertex digraph $D$.
An arc $ (v_{i},v_{j})\in E(D) $ is a \textit{backward arc of $D$ under} $ \theta $ if $ i > j $. We say that a vertex $ v_{j} $ is \textit{between} two vertices $ v_{i},v_{k} $ under $ \theta = (v_{1},...,v_{n}) $ if $ i < j < k $ or $ k < j < i $. 
 The graph of backward arcs under $ \theta $, denoted by $ B(D,\theta) $, is the undirected graph that has vertex set $V(D)$, and $ v_{i}v_{j} \in E(B(D,\theta)) $ if and only if $ (v_{i},v_{j}) $ or $ (v_{j},v_{i}) $ is a backward arc of $D$ under $ \theta $. A tournament $S$ on $p$ vertices with $V(S)= \lbrace u_{1},u_{2},...,u_{p}\rbrace$ is a \textit{right star} (resp. \textit{left star}) (resp. \textit{middle star}) if there exist an ordering $\theta^{*} = (u_{1},u_{2},...,u_{p})$ of its vertices such that the backward arcs of $S$ under $\theta^{*}$ are $(u_{p},u_{i})$ for $i=1,...,p-1$ (resp. $(u_{i},u_{1})$ for  $i=2,...,p$) (resp. $(u_{i},u_{r})$ for $i= r+1,...,p$ and $(u_{r},u_{i})$ for $i=1,...,r-1$, where $2\leq r\leq p-1$). In this case we write $S = \lbrace u_{1},u_{2},...,u_{p}\rbrace$ and we call $\theta^{*} = (u_{1},u_{2},...,u_{p})$ a \textit{right star ordering} (resp. \textit{left star ordering}) (resp. \textit{middle star ordering}) of $S$, $u_{p}$ (resp. $u_{1}$) (resp. $u_{r}$) the \textit{center of} $S$, and $u_{1},...,u_{p-1}$ (resp. $u_{2},...,u_{p}$) (resp. $u_{1},...,u_{r-1},u_{r+1},...,u_{p}$) the \textit{leaves of} $S$. A \textit{star} is a left star or a right star or a middle star. A \textit{star ordering} is a left star ordering or a right star ordering or a middle star ordering. Note that in the case $p=2$ we may choose arbitrarily any one of the two vertices to be the center of the star, and the other vertex is then considered to be the leaf. A \textit{frontier star} is a left star or a right star (note that a frontier star is not a middle star, a frontier star is either left or right).  
   A \textit{star} $S=\lbrace v_{i_{1}},...,v_{i_{t}}\rbrace$ \textit{of $D$ under $\theta$} (where $i_{1}<...<i_{t}$) is the subdigraph of $D$ induced by $\lbrace v_{i_{1}},...,v_{i_{t}}\rbrace$ such that $S$ is a star and $S$ has the star ordering $ (v_{i_{1}},...,v_{i_{t}})$ under $\theta$ (i.e $(v_{i_{1}},...,v_{i_{t}})$ is the restriction of $\theta$ to $V(S)$ and $ (v_{i_{1}},...,v_{i_{t}})$ is a star ordering of $S$).

A tournament $T$ is a \textit{galaxy} if there exists an ordering $\theta$ of its vertices such that $V(T)$ is the disjoint union of $V(Q_{1}),...,V(Q_{l}),X$ where $Q_{1},...,Q_{l}$ are the frontier stars of $T$ under $\theta$, and for every $x\in X$, $\lbrace x \rbrace$ is a singleton component of $B(T,\theta)$, and no center of a star is between leaves of another star under $\theta$. In this case we also say that $T$ \textit{is a \textit{galaxy} under $\theta$}. If $X=\phi$, we say that $T$ is a \textit{regular galaxy under $\theta$}. \\
The condition that no center of a star appears in the ordering between leaves of another star is necessary to make the proof of the following theorem work: 
\begin{theorem}\cite{polll}
Every galaxy satisfies the Erd\"{o}s-Hajnal conjecture.
\end{theorem}
It is not known whether the conjecture is still satisfied  if the condition concerning the centers of the stars in galaxies is abandoned. A tournament $T$ is a \textit{nebula} if there exists an ordering $\theta$ of its vertices such that $V(T)$ is the disjoint union of $V(Q_{1}),...,V(Q_{l}),X$ where $Q_{1},...,Q_{l}$ are the stars of $T$ under $\theta$, and for every $x\in X$, $\lbrace x \rbrace$ is a singleton component of $B(T,\theta)$ (the star does not have to be necessarily left or right, and there is no condition concerning the location of the centers of the stars). We also say that $T$ \textit{is a \textit{nabula} under $\theta$}. If $X=\phi$, we say that $T$ is a \textit{regular nebula}. Notice that every galaxy is a nebula. The following is still a conjecture:
\begin{conjecture} \label{j}
Every nebula satisfies the Erd\"{o}s-Hajnal conjecture.
\end{conjecture}
In $2015$ Choromanski proved Conjecture \ref{j} for every \textit{constellation} \cite{kg} (every constellation is a nebula). In the constellation case we do not have middle stars and we have center of a star between leaves of another star under some conditions. Constellations are fully characterized in \cite{kg}. Recently we prove Conjecture \ref{j} for every \textit{galaxy with spiders} \cite{sss}. In galaxies with spiders we have middle stars also and we have center of a star between leaves of another star under some conditions. Galaxies with spiders are fully caracterized in \cite{sss} (every galaxy with spiders is a nebula).\vspace{1.5mm}\\
In what follows we define formally the family of \textit{super nebulas} and the family of \textit{super triangular galaxies}.\\
In order to define super nebulas we need first to define \textit{$2$-nebulas} and \textit{super $2$-nebulas}. Let $s$ be a $\lbrace 0,1 \rbrace$$-$vector. Denote $s_{c}$ the vector obtained from $s$ by replacing every subsequence of consecutive $1'$s by single $1$. Let  $T$ be a regular nebula under $\theta = (v_{1},...,v_{n})$ and let $Q_{1},...,Q_{l}$ be the stars of $T$ under $\theta$. Let $s^{T,\theta}$ be the $\lbrace 0,1\rbrace$$-$vector such that $s^{T,\theta}_{i}=1$ if and only if $v_{i}$ is a leaf of one of the stars of $T$ under $\theta$. A tournament $T$ is a \textit{$2$-nebula} if it is a regular nebula under $\theta$ and besides $V(T)$ is the disjoint union of $V(Q_{1})$ and $V(Q_{2})$, where $Q_{1}$ and $Q_{2}$ are the frontier stars of $T$ under $\theta$, and $s^{T,\theta}_{c}$ is one of the following vectors: $ (0,0,1)$, or $(0,1,0)$, or $ (1,0,0)$. A $2$-nebula $T$ under $\theta$ is a \textit{left} (resp. \textit{middle}) (resp. \textit{right}) $2$\textit{-nebula under} $\theta$  if $s^{T,\theta}_{c} =(1,0,0)$ (resp. $s^{T,\theta}_{c} =(0,1,0)$) (resp. $s^{T,\theta}_{c} =(0,0,1)$). A tournament $H$ is a \textit{super left $2$-nebula} (resp. \textit{super middle $2$-nebula}) (resp. \textit{super right $2$-nebula}) under $\theta = (v_{1},...,v_{n})$ if it is obtained from a left $2$-nebula (resp. middle $2$-nebula) (resp. right $2$-nebula) $T$ under $\theta$ by reversing the orientation of the arc joining the centers of the two frontier stars of $T$ under $\theta$. In this case we write $H =\lbrace v_{1},...,v_{n}\rbrace$, we call $\theta$ a \textit{super left $2$-nebula ordering} (resp. \textit{super middle $2$-nebula ordering}) (resp. \textit{super right $2$-nebula ordering}) \textit{of $H$}, the leaves of the two frontier stars of $T$ under $\theta$ are called the \textit{leaves of} $H$, and the centers of the two frontier stars of $T$ under $\theta$ are called the \textit{centers of} $H$. A \textit{super $2$-nubula} is a super left $2$-nebula or a super middle $2$-nebula or a super right $2$-nebula. A \textit{super $2$-nubula ordering} is a super left $2$-nebula ordering or a super middle $2$-nebula ordering or a super right $2$-nebula ordering. Let $v_{q_{1}}$ and $v_{q_{2}}$ be the centers of $Q_{1}$ and $Q_{2}$ respectively. The leaves of $Q_{1}$ (resp. $Q_{2}$) are called the \textit{leaves of $H$ incident to the center $v_{q_{1}}$ (resp. $v_{q_{2}}$) of $H$}.

\begin{figure}[h]
	\centering
	\includegraphics[width=0.5\linewidth]{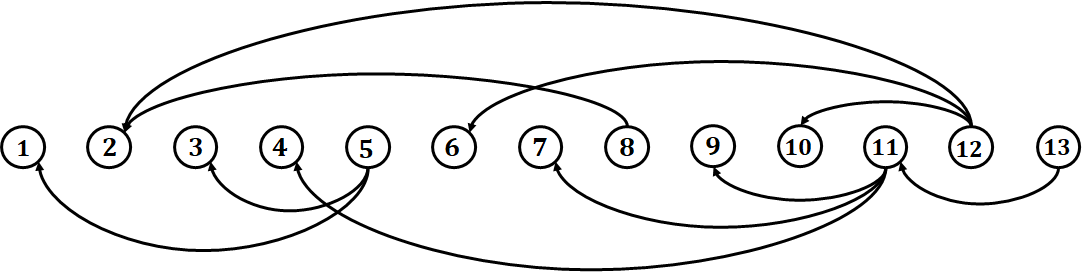}
	\caption{Super nebula consisting of one super $2$-nebula and two stars. It is drawn under its super nebula ordering. All the arcs that are not drawn are forward.}
	\label{fig:supernebula}
\end{figure}

 Let $\theta = (v_{1},...,v_{n})$ be an ordering of the vertex set $V(G)$ of an $n$$-$vertex tournament $G$. A \textit{super $2$-nebula $T =\lbrace v_{i_{1}},...,v_{i_{t}}\rbrace$ of $G$ under} $\theta$ is the subtournament of $G$ induced by $\lbrace v_{i_{1}},...,v_{i_{t}}\rbrace$ such that $T$ has the super $2$-nebula ordering $ (v_{i_{1}},...,v_{i_{t}})$ under $\theta$ (note that $i_{1}<...<i_{t}$).
Now we are ready to define super nebulas:\vspace{1.5mm}

A tournament $T$ is a \textit{super nebula} if there exists an ordering $\theta$ of its vertices such that $V(T)$ is the disjoint union of $V(Q_{1}),...,V(Q_{m}),V(\Sigma_{1}),...,V(\Sigma_{l}),X$ where $Q_{1},...,Q_{m}$ are the stars of $T$ under $\theta$, $\Sigma_{1},...,\Sigma_{l}$ are the super $2$-nebulas of $T$ under $\theta$, no center of a star is between leaves of a super $2$-nebula under $\theta$, no center of a super $2$-nebula is between leaves of another super $2$-nebula under $\theta$, and for every $x\in X$, $\lbrace x \rbrace$ is a singleton component of $B(T,\theta)$. In this case we say that $T$ is a \textit{super nebula under} $\theta$ and $\theta$ is called a \textit{super nebula ordering of $T$}. If $X=\phi$ then $T$ is called \textit{regular super nebula} (see Figure \ref{fig:supernebula}). Obviously, every nebula is a super nebula. 

\begin{figure}[h]
	\centering
	\includegraphics[width=0.30\linewidth]{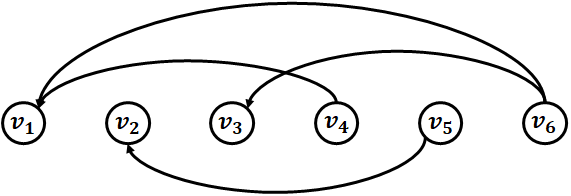}
	\caption{Tournament $K_{6}$ drawn under its canonical ordering. All the arcs that are not drawn are forward.}
	\label{fig:K6tour}
\end{figure}

Denote by $K_{6}$ the six-vertex tournament with $V(K_{6}) = \lbrace v_{1},v_{2},v_{3},v_{4},v_{5},v_{6} \rbrace$ such that under ordering $(v_{1},v_{2},v_{3},v_{4},v_{5},v_{6})$ of its vertices the set of backward arcs is: $ \lbrace (v_{4},v_{1}),(v_{6},v_{3}),(v_{6},v_{1}),(v_{5},v_{2}) \rbrace$. We call this ordering of vertices of $K_{6}$ the \textit{canonical ordering of} $K_{6}$. $K_{6}$ is the only tournament on at most six vertices for which the conjecture is still open \cite{bnmm}. Notice that $K_{6}$ is obviously a super nebula and its canonical ordering is its super nebula ordering (see Figure \ref{fig:K6tour}).
\\Let $\theta = (v_{1},...,v_{n})$ be an ordering of the vertex set $V(T)$ of an $n$$-$vertex tournament $T$. A \textit{triangle} $C =\lbrace v_{i_{1}},v_{i_{2}},v_{i_{3}}\rbrace$ \textit{of $T$ under $\theta$} is a transitive subtournament of $T$ induced by $\lbrace v_{i_{1}},v_{i_{2}},v_{i_{3}}\rbrace$ such that $(v_{i_{3}},v_{i_{2}},v_{i_{1}})$ is its transitive ordering (i.e $v_{i_{1}}\leftarrow v_{i_{2}}$ and $\lbrace v_{i_{1}},v_{i_{2}}\rbrace\leftarrow v_{i_{3}}$), and $i_{1}<i_{2}<i_{3}$.
 We call $ v_{i_{2}} $ the \textit{center of the $C$}, $ v_{i_{1}}$ the \textit{ left exterior of $C$}, and $ v_{i_{3}} $ the \textit{right exterior of $C$}.\vspace{1.5mm}
 
\begin{figure}[h]
	\centering
	\includegraphics[width=0.35\linewidth]{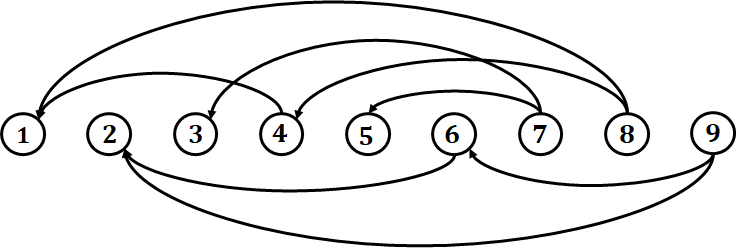}
	\caption{Central triangular galaxy consisting of two triangles and one star. All the non drawn arcs are forward.}
	\label{fig:triangle}
\end{figure}

A tournament $T$ is a \textit{triangular galaxy}  if there exist an ordering $\theta$ of its vertices such that $V(T)$ is the disjoint union of $V(\Delta_{1}),...,V(\Delta_{l}),X$ where $\Delta_{1},...,\Delta_{l}$ are the triangles of $T$ under $\theta$, and $T$$\mid$$X$ is a galaxy under $\theta_{X}$ where $\theta_{X}$ is the restriction of $\theta$ to $X$, and no  vertex of a triangle is between leaves of a star of $T$ under $\theta$. In this case we say that $T$ is \textit{triangular galaxy under $\theta$} and $\theta$ is a \textit{triangular galaxy ordering of $T$}. If $T$$\mid$$X$ is a regular galaxy under $\theta_{X}$ then we say that $T$ is a \textit{regular triangular galaxy under $\theta$}. If $l=1$ then $T$ is called \textit{$\Delta$galaxy under} $\theta$ and  $\theta$ is a \textit{$\Delta$galaxy ordering of $T$}. If $l=1$ and $T$$\mid$$X$ is a regular galaxy under $\theta_{X}$ then $T$ is called \textit{regular $\Delta$galaxy under} $\theta$. If for every $x\in X$, $\lbrace x \rbrace$ is a singleton component of $B(T,\theta)$ we say that $T$ is a \textit{triangular tournament}. If $X=\phi$, we say that $T$ is a \textit{regular triangular tournament}. If the condition concerning the location of the vertices of the triangles is weakend such that the centers (resp. right exteriors) (resp. left exteriors) of the triangles are allowed to be  between leaves of a star of $T$ under $\theta$ then $T$ is called  \textit{central triangular galaxy under $\theta$} (resp. \textit{right triangular galaxy under $\theta$}) (resp. \textit{left triangular galaxy under $\theta$}) and  $\theta$ is called a \textit{central triangular galaxy ordering of $T$} (resp. \textit{right triangular galaxy ordering of $T$}) (resp. \textit{left triangular galaxy ordering of $T$}) (see Figure\ref{fig:triangle}). If the condition concerning the vertices of the triangles is abandoned or weakened then $T$ is called a \textit{super triangular galaxy under $\theta$}.
\vspace{2mm}\\
Unfortunately the following are still open:
\begin{problem}
Does the tournament $K_{6}$ satisfy EHC?
\end{problem}
\begin{problem}
Does every nebula (resp. super nebula) satisfy EHC?
\end{problem}
\begin{problem}
Does every triangular tournament (resp. triangular galaxy) (resp. $\Delta$galaxy) (resp. super triangular galaxy) satisfy EHC?
\end{problem}
However if we exclude both:
\begin{itemize}
\item An arbitrary super nebula and an arbitrary $\Delta$galaxy, or
\item An arbitrary central triangular galaxy and $K_{6}$, or
\item An arbitrary middle $\Sigma$-galaxy and an arbitrary central triangular galaxy, or
\item An arbitrary left $\Sigma$-galaxy and an arbitrary  left triangular galaxy, or
\item An arbitrary right $\Sigma$-galaxy and an arbitrary  right triangular galaxy
\end{itemize}
then the we prove that the conjecture is satisfied. More results in this setting concerning nebulas, super nebulas, and super triangular galaxies are formally stated in Section \ref{k}.\vspace{2mm}\\
Middle (resp. left) (resp. right) $\Sigma$-galaxies are defined in Section \ref{k}. Every middle (resp. left) (resp. right) $\Sigma$-galaxy is a super nebula.\\
In \cite{ssss} we prove that some $7$$-$vertex super triangular galaxy tournaments satisfy $EHC$ and in \cite{sss} we prove that if the super nebula tournament is an \textit{asterism} then it satisfies $EHC$. \vspace{2mm}\\
The main results of this paper are the following:
\begin{theorem} \label{p}
If $\mathcal{N}$ is a super nebula and $\mathcal{G}$ is a $\Delta$galaxy, then  $\lbrace \mathcal{N},\mathcal{G}\rbrace$ satisfy the Erd\"{o}s-Hajnal Conjecture.  
\end{theorem} 
\begin{theorem}\label{t}
If $H$ is a central triangular galaxy, then $\lbrace K_{6},H\rbrace$ satisfy the Erd\"{o}s-Hajnal Conjecture.  
\end{theorem}
The results concerning the extensions of Theorem \ref{p} and Theorem \ref{t} are formally stated in Section \ref{k}.\vspace{2mm}\\ 
Let $T$ be a nebula under $\theta$ and let $Q_{1},...,Q_{l}$ be the stars of $T$ under $\theta$. $T$ is a \textit{left nebula} (resp. \textit{right nebula}) (resp. \textit{central nebula}) under $\theta$ if for all $i\in \lbrace 1,...,l\rbrace$, $Q_{i}$ is a $3$-vertex left star (resp. $3$-vertex right star) (resp. $3$-vertex middle star). In the setting of forbidding pairs of tournaments, Choromanski proved \cite{kgg} the following: 
\begin{theorem}\cite{kgg}
If $H_{1}$ and $H_{2}$ are: a left nebula and a  right nebula, or: a left nebula and a  central nebula, or: a right nebula and a  central nebula, then $\lbrace H_{1},H_{2}\rbrace$ satisfies the Erd\"{o}s-Hajnal Conjecture. 
\end{theorem} 
 This paper is organized as follows:\vspace{1mm}\\
$\bullet$ In section $2$ we give some definitions and preliminary lemmas, and we prove some lemmas needed in the proof of the main results in this paper.\\
$\bullet$ In section $3$ we introduce some tools useful in the proof of Theorem \ref{p} and we prove Theorem \ref{p}.\\
$\bullet$ In section $4$ we introduce some definitions and we prove Theorem \ref{t}.\\
$\bullet$ In section $5$ we give extensions of our results.
\section{Definitions and Preliminary Lemmas}
Denote by $tr(T)$ the largest size of a transitive subtournament of a tournament $T$. For $X \subseteq V(T)$, write $tr(X)$ for $tr(T$$\mid$$X)$. Let $X, Y \subseteq V(T)$ be disjoint. Denote by $e_{X,Y}$ the number of directed arcs $(x,y)$, where $x \in X$ and $y \in Y$. The \textit{directed density from $X$ to} $Y$ is defined as $d(X,Y) = \frac{e_{X,Y}}{\mid X \mid.\mid Y \mid} $. We call $T$ $ \epsilon - $\textit{critical} for $ \epsilon > 0 $ if $tr(T) < $ $ \mid $$T$$ \mid^{\epsilon} $ but for every proper subtournament $S$ of $T$ we have: $tr(S) \geq $ $ \mid $$S$$ \mid^{\epsilon}. $
\begin{lemma} \cite{ml} \label{h}
Every tournament on $2^{k-1}$ vertices contains a transitive subtournament of size at least $k$.
 \end{lemma}
\begin{lemma} \cite{polll} \label{e} For every $N$ $ > 0 $, there exists $ \epsilon(N) > 0 $ such that for every $ 0 < \epsilon < \epsilon(N)$ every $ \epsilon $$-$critical tournament $T$ satisfies $ \mid $$T$$\mid$ $ \geq N$.
\end{lemma} 
\begin{lemma} \cite{polll} \label{v} Let $T$ be an $ \epsilon $$-$critical tournament with $\mid$$T$$\mid = n$ and $\epsilon$,$c > 0 $ be constants such that $ \epsilon < log_{\frac{c}{2}}(\frac{1}{2}). $ Then for every two disjoint subsets $X, Y \subseteq V(T)$ with $ \mid $$X$$ \mid$ $ \geq cn, \mid $$Y$$ \mid$ $ \geq cn $ there exist an integer $k \geq \frac{cn}{2} $ and vertices $ x_{1},...,x_{k} \in X $ and $ y_{1},...,y_{k} \in Y $ such that $ y_{i} $ is adjacent to $ x_{i} $ for $i = 1,...,k. $
\end{lemma}
\begin{lemma} \cite{polll} \label{f} Let $T$ be an $ \epsilon $$-$critical tournament with $\mid$$T$$\mid$ $ =n$ and $\epsilon$,$c,f > 0 $ be constants such that  $ \epsilon < log_{c}(1 - f)$. Then for every $A \subseteq V(T)$ with $ \mid $$A$$ \mid$ $ \geq cn$ and every transitive subtournament $G$ of $T$ with $ \mid $$G$$ \mid$ $\geq f.tr(T)$ and $V(G) \cap A = \phi$, we have: $A$ is not complete from $V(G)$ and $A$ is not complete to $V(G)$.
\end{lemma}
\begin{lemma}\label{r}
Let $f_{1},...,f_{m},l_{1},...,l_{t},c,\epsilon > 0$ be constants, where $0 <  f_{1},...,f_{m},l_{1},...,l_{t},c < 1$ and $0 < \epsilon < min\lbrace log_{\frac{c}{2(m+t)}}(1-f_{1}),..., log_{\frac{c}{2(m+t)}}(1-f_{m}),log_{\frac{c}{2(m+t)}}(1-l_{1}),..., log_{\frac{c}{2(m+t)}}(1-l_{t})\rbrace$. Let $T$ be an $ \epsilon $$-$critical tournament with $\mid$$T$$\mid$ $ =n$, and let $S_{1},...,S_{m},P_{1},...,P_{t}$ be $m+t$ disjoint transitive subtournaments of $T$ with $m,t\in\mathbb{N}$, $ \mid $$S_{i}$$ \mid$ $\geq f_{i}.tr(T)$ for $i=1,...,m$, and $ \mid $$P_{i}$$ \mid$ $\geq l_{i}.tr(T)$ for $i=1,...,t$. Let $A \subseteq V(T) \backslash ((\bigcup_{i=1}^{m} V(S_{i}))\cup (\bigcup_{i=1}^{t} V(P_{i})))$ with $ \mid $$A$$ \mid$ $ \geq cn$. Then there exist vertices $s_{1},...,s_{m},p_{1},...,p_{t},g$ such that $g\in A$, $s_{i}\in S_{i}$ for $i=1,...,m$, $p_{i}\in P_{i}$ for $i=1,...,t$, and $\lbrace s_{1},...,s_{m} \rbrace\leftarrow g \leftarrow\lbrace p_{1},...,p_{t} \rbrace$.  
\end{lemma}
\begin{proof}
Let $A_{i} \subseteq A$ such that $A_{i}$ is complete from $S_{i}$ for $i = 1,...,m$ and let $A^{i} \subseteq A$ such that $A^{i}$ is complete to $P_{i}$ for $i = 1,...,t$. Let $1\leq j \leq m$. If $\mid$$A_{j}$$\mid \geq \frac{\mid A \mid}{2(m+t)}\geq \frac{c}{2(m+t)}n$, then this will contradicts Lemma \ref{f} since $\mid$$S_{j}$$\mid \geq f_{j}tr(T)$ and $\epsilon < log_{\frac{c}{2(m+t)}}(1-f_{j})$. Then $\forall i \in \lbrace 1,...,m \rbrace$, $\mid$$A_{i}$$\mid < \frac{\mid A \mid}{2(m+t)}$. Similarly we prove that $\forall i \in \lbrace 1,...,t \rbrace$, $\mid$$A^{i}$$\mid < \frac{\mid A \mid}{2(m+t)}$. Let $A^{*} = A\backslash ((\bigcup_{i=1}^{m}A_{i})\cup (\bigcup_{i=1}^{t}A^{i}))$, then $\mid$$A^{*}$$\mid > \mid $$A$$ \mid-(m+t).\frac{\mid A \mid}{2(m+t)} \geq \frac{\mid A \mid}{2}$. Then $A^{*} \neq \phi$. Fix $g\in A^{*}$. So there exist vertices $s_{1},...,s_{m},p_{1},...,p_{t}$ such that $s_{i}\in S_{i}$ for $i=1,...,m$, $p_{i}\in P_{i}$ for $i=1,...,t$ and $\lbrace s_{1},...,s_{m} \rbrace\leftarrow g \leftarrow\lbrace p_{1},...,p_{t} \rbrace$.   $\blacksquare$ 
\end{proof}
\begin{lemma} \label{s}
Let $f,c,t,\epsilon > 0$ be constants, where $0 <  f,c < 1$, $t$ a positive integer, and $0 < \epsilon < min\lbrace log_{\frac{c}{2t}}(1-f), log_{\frac{c}{4}}(\frac{1}{2})\rbrace$. Let $T$ be an $ \epsilon $$-$critical tournament with $\mid$$T$$\mid$ $ =n$, and   let $S_{1},...,S_{m},S_{m+1},...,S_{t}$ be $t$ disjoint transitive subtournaments of $T$ with $ \mid $$S_{i}$$ \mid$ $\geq f.tr(T)$ for $i=1,...,t$. Let $A_{1},A_{2}$ be two disjoint subsets of $V(T)$ with $ \mid $$A_{1}$$ \mid$ $\geq cn$, $ \mid $$A_{2}$$ \mid$ $\geq cn$, and $A_{1},A_{2} \subseteq V(T) \backslash (\bigcup_{i=1}^{t}S_{i})$. Then there exist vertices $x,y,s_{1},...,s_{t}$ such that $x\in A_{1}, y\in A_{2}, s_{i}\in S_{i}$ for $i=1,...,t$, $x\leftarrow\lbrace y,s_{1},...,s_{m}\rbrace$, and $\lbrace s_{m+1},...,s_{t}\rbrace\leftarrow y$. Similarly there exist vertices $p,q,u_{1},...,u_{t}$ such that $p\in A_{1}, q\in A_{2}, u_{i}\in S_{i}$ for $i=1,...,t$, $\lbrace u_{1},...,u_{m},p\rbrace\leftarrow q$, and $\lbrace u_{m+1},...,u_{t}\rbrace\leftarrow p$. Similarly there exist vertices $g,v,z_{1},...,z_{t}$ such that $g\in A_{1}, v\in A_{2}, z_{i}\in S_{i}$ for $i=1,...,t$, $g\leftarrow\lbrace v,z_{1},...,z_{m}\rbrace$, and $v\leftarrow\lbrace z_{m+1},...,z_{t}\rbrace$.
\end{lemma}
\begin{proof}
We will prove only the first statement  because the rest can be proved analogously. Let $A_{1}^{*} = \lbrace x \in A_{1}; \exists s_{i} \in S_{i}$ for $i=1,...,m$ and $x \leftarrow \lbrace s_{1},...,s_{m}\rbrace \rbrace$ and let $A_{2}^{*} = \lbrace y \in A_{2}; \exists s_{i} \in S_{i}$ for $i=m+1,...,t$ and $\lbrace s_{m+1},...,s_{t}\rbrace \leftarrow y \rbrace$. Since $ \epsilon < log_{\frac{c}{2t}}(1-f)$, then $\mid$$A_{1}^{*}$$\mid$ $> \frac{\mid A_{1} \mid}{2} \geq \frac{c}{2}n$ and $\mid$$A_{2}^{*}$$\mid$ $> \frac{\mid A_{2} \mid}{2} \geq \frac{c}{2}n$. Now since $\epsilon < log_{\frac{c}{4}}(\frac{1}{2})$, then Lemma \ref{v} implies that $\exists k \geq \frac{c}{4}n$, $\exists x_{1},...,x_{k} \in A_{1}^{*}$, $\exists y_{1},...,y_{k} \in A_{2}^{*}$, such that $x_{i} \leftarrow y_{i}$ for $i = 1,...,k$. Fix $i_{1}\in \lbrace 1,...,k\rbrace$.  So, there exist vertices $s_{1},...,s_{t}$ such that $ s_{i}\in S_{i}$ for $i=1,...,t$, $x_{i_{1}}\leftarrow\lbrace y_{i_{1}},s_{1},...,s_{m}\rbrace$, and $\lbrace s_{m+1},...,s_{t}\rbrace\leftarrow y_{i_{1}}$. $\blacksquare$ 
\end{proof}
\begin{lemma} \cite{polll} \label{b} Let $A_{1},A_{2}$ be two disjoint sets such that $d(A_{1},A_{2}) \geq 1-\lambda$ and let $0 < \eta_{1},\eta_{2} \leq 1$. Let $\widehat{\lambda} = \frac{\lambda}{\eta_{1}\eta_{2}}$. Let $X \subseteq A_{1}, Y \subseteq A_{2}$ be such that $\mid$$X$$\mid$ $\geq \eta_{1} \mid$$A_{1}$$\mid$ and $\mid$$Y$$\mid$ $\geq \eta_{2} \mid$$A_{2}$$\mid$. Then $d(X,Y) \geq 1-\widehat{\lambda}$. 
\end{lemma}
The following is introduced in \cite{bnmm}.\\
Let $ c > 0, 0 < \lambda < 1 $ be constants, and let $w$ be a $ \lbrace 0,1 \rbrace - $ vector of length $ \mid $$w$$ \mid $. Let $T$ be a tournament with $ \mid $$T$$ \mid$ $ = n. $ A sequence of disjoint subsets $ \chi = (S_{1}, S_{2},..., S_{\mid w \mid}) $ of $V(T)$ is a smooth $ (c,\lambda, w)- $structure if:\\
$\bullet$ whenever $ w_{i} = 0 $ we have $ \mid $$S_{i}$$ \mid$ $ \geq cn $ (we say that $ S_{i} $ is a \textit{linear set}).\\
$\bullet$ whenever $ w_{i} = 1 $ the tournament $T$$\mid$$ S_{i} $ is transitive and $ \mid $$S_{i}$$ \mid$ $ \geq c.tr(T) $ (we say that $ S_{i} $ is a \textit{transitive set}).\\
$\bullet$ $ d(\lbrace v \rbrace, S_{j}) \geq 1 - \lambda $ for $v \in S_{i} $ and $ d(S_{i}, \lbrace v \rbrace) \geq 1 - \lambda $ for $v \in S_{j}, i < j $ (we say that $\chi$ is \textit{smooth}).
\begin{theorem} \cite{bnmm} \label{i}
Let $S$ be a tournament, let $w$ be a $ \lbrace 0,1 \rbrace - $vector, and let $ 0 < \lambda_{0} < \frac{1}{2} $ be a constant. Then there exist $ \epsilon_{0}, c_{0} > 0 $ such that for every $ 0 < \epsilon < \epsilon_{0} $, every $ S- $free $ \epsilon $$- $critical tournament contains a smooth $ (c_{0}, \lambda_{0},w)- $structure.
\end{theorem} 

Let $(S_{1},...,S_{\mid w \mid})$ be a smooth $(c,\lambda ,w)$$-$structure of a tournament $T$, let $i \in \lbrace 1,...,\mid$$w$$\mid \rbrace$, and let $v \in S_{i}$. For $j\in \lbrace 1,2,...,\mid$$w$$\mid \rbrace \backslash \lbrace i \rbrace$, denote by $S_{j,v}$ the set of the vertices of $S_{j}$ adjacent from $v$ for $j > i$ and adjacent to $v$ for $j<i$.
\begin{lemma} \label{g} Let $0<\lambda<1$, $0<\gamma \leq 1$ be constants and let $w$ be a $\lbrace 0,1 \rbrace$$-$vector. Let $(S_{1},...,S_{\mid w \mid})$ be a smooth $(c,\lambda ,w)$$-$structure of a tournament $T$ for some $c>0$. Let $j\in \lbrace 1,...,\mid$$w$$\mid \rbrace$. Let $S_{j}^{*}\subseteq S_{j}$ such that $\mid$$S_{j}^{*}$$\mid$ $\geq \gamma \mid$$S_{j}$$\mid$ and let $A= \lbrace x_{1},...,x_{k} \rbrace \subseteq \displaystyle{\bigcup_{i\neq j}S_{i}}$ for some positive integer $k$. Then $\mid$$\displaystyle{\bigcap_{x\in A}S^{*}_{j,x}}$$\mid$ $\geq (1-k\frac{\lambda}{\gamma})\mid$$S_{j}^{*}$$\mid$. In particular $\mid$$\bigcap_{x\in A}S_{j,x}$$\mid$ $\geq (1-k\lambda)\mid$$S_{j}$$\mid$.
\end{lemma}
\begin{proof}
The proof is by induction on $k$. without loss of generality assume that $x_{1} \in S_{i}$ and $j<i$.\\ Since $\mid$$S_{j}^{*}$$\mid$ $\geq \gamma \mid$$S_{j}$$\mid$ then by Lemma \ref{b}, $d(S^{*}_{j},\lbrace x_{1}\rbrace) \geq 1-\frac{\lambda}{\gamma}$. So $1-\frac{\lambda}{\gamma} \leq d(S^{*}_{j},\lbrace x_{1}\rbrace) = \frac{\mid S^{*}_{j,x_{1}}\mid}{\mid S_{j}^{*}\mid}$.\\ Then $\mid$$S^{*}_{j,x_{1}}$$\mid$ $\geq (1-\frac{\lambda}{\gamma})$$\mid$$S_{j}^{*}$$\mid$ and so true for $k=1$.
Suppose the statement is true for $k-1$.\\ $\mid$$\displaystyle{\bigcap_{x\in A}S^{*}_{j,x}}$$\mid$ $=\mid$$(\displaystyle{\bigcap_{x\in A\backslash \lbrace x_{1}\rbrace}S^{*}_{j,x}})\cap S^{*}_{j,x_{1}}$$\mid$ $= \mid$$\displaystyle{\bigcap_{x\in A\backslash \lbrace x_{1}\rbrace}S^{*}_{j,x}}$$\mid$ $+$ $\mid$$S^{*}_{j,x_{1}}$$\mid$ $- \mid$$(\displaystyle{\bigcap_{x\in A\backslash \lbrace x_{1}\rbrace}S^{*}_{j,x}})\cup S^{*}_{j,x_{1}}$$\mid$ $\geq (1-(k-1)\frac{\lambda}{\gamma})\mid$$S_{j}^{*}$$\mid$ $+$ $(1-\frac{\lambda}{\gamma})\mid$$S_{j}^{*}$$\mid$ $-$ $\mid$$S_{j}^{*}$$\mid$ $= (1-k\frac{\lambda}{\gamma})\mid$$S_{j}^{*}$$\mid$.\hspace{2mm} $\blacksquare$       
\end{proof}
\section{Super Nebulas and $\Delta$galaxies}
\subsection{Definitions and Tools} \label{c}

In this section we introduce the notion of \textit{bad triplets}, the notion of \textit{key tournaments}, and the notion of \textit{corresponding structures under a tournament and an ordering of its vertices} that will be very crucial in our later analysis.\vspace{2mm}

Let $s$ be a $\lbrace 0,1 \rbrace$$-$vector. Denote $s_{c}$ the vector obtained from $s$ by replacing every subsequence of consecutive $1'$s by single $1$. 
Let $\mathcal{N}$ be a regular super nebula under $\theta = (v_{1},...,v_{n})$. Let $Q_{1},...,Q_{m}$ be the stars of $\mathcal{N}$ under $\theta $ and let $\Sigma_{1},...,\Sigma_{l}$ be the super $2$-nebulas of $\mathcal{N}$ under $\theta$. Let $s^{\mathcal{N},\theta}$ be the $\lbrace 0,1\rbrace$$-$vector such that $s^{\mathcal{N},\theta}_{i}=1$ if and only if $v_{i}$ is a leaf of one of the stars of $\mathcal{N}$ under $\theta$ or a leaf of one of the super $2$-nebulas of $\mathcal{N}$ under $\theta$. Let $\omega =s^{\mathcal{N},\theta}_{c}$ and let $i_{r}$ be such that $\omega_{i_{r}}=1$. Let $j$ be such that $s^{\mathcal{N},\theta}_{j}=1$. We say that $s^{\mathcal{N},\theta}_{j}$ \textit{corresponds to} $\omega_{i_{r}}$ if $s^{\mathcal{N},\theta}_{j}$ belongs to the subsequence of consecutive $1'$s that is replaced by the entry $\omega_{i_{r}}$. Define $\mathcal{R}_{i_{r}} = \lbrace v_{i}\in V(\mathcal{N}); s^{\mathcal{N},\theta}_{i}=1$ and $s^{\mathcal{N},\theta}_{i}$ corresponds to $\omega_{i_{r}}\rbrace$. Let $i_{1},...,i_{t}$ be the non zero entries of $\omega$. Define $\mathcal{R}=\bigcup_{j=1}^{t}\mathcal{R}_{i_{j}}$. Notice that $\mathcal{R}$ is the set of the leaves of the stars and super $2$-nebulas of $\mathcal{N}$ under $\theta$. Let $1\leq i\leq m$. Let $L_{i}$ be the set of leaves of the star $Q_{i}$ and let $v_{q_{i}}$ be the center of $Q_{i}$. For all $v_{j},v_{k}\in L_{i}$, if $v_{j}\in \mathcal{R}_{i_{r}}$ and $v_{k}\notin \mathcal{R}_{i_{r}}$ for some $1\leq r \leq t$ then the triplet $\lbrace v_{q_{i}},v_{j},v_{k}\rbrace$ is called a \textit{bad triplet of $Q_{i}$}. Denote by $B_{i}$ the set of all bad triplets of $Q_{i}$. Define $\mathcal{T}^{\mathcal{N},\theta}=\lbrace \mathcal{T}_{j}=\lbrace v_{j_{1}},v_{j_{2}},v_{j_{3}}\rbrace \subseteq \bigcup_{i=1}^{l}V(\Sigma_{i}); \mid$$E(B(\mathcal{N}$$\mid$$\lbrace v_{j_{1}},v_{j_{2}},v_{j_{3}}\rbrace$$,\theta_{\mathcal{T}_{j}}))$$\mid =2$ where $\theta_{\mathcal{T}_{j}}$ is the restriction of $\theta$ to $\mathcal{T}_{j}\rbrace$. Denote by $B^{\mathcal{N},\theta}=(\bigcup_{i=1}^{m}B_{i})\cup \mathcal{T}^{\mathcal{N},\theta}$ the set of all \textit{bad triplets of $\mathcal{N}$ under $\theta$}.\vspace{3.5mm}

\begin{figure}[h]
	\centering
	\includegraphics[width=0.6\linewidth]{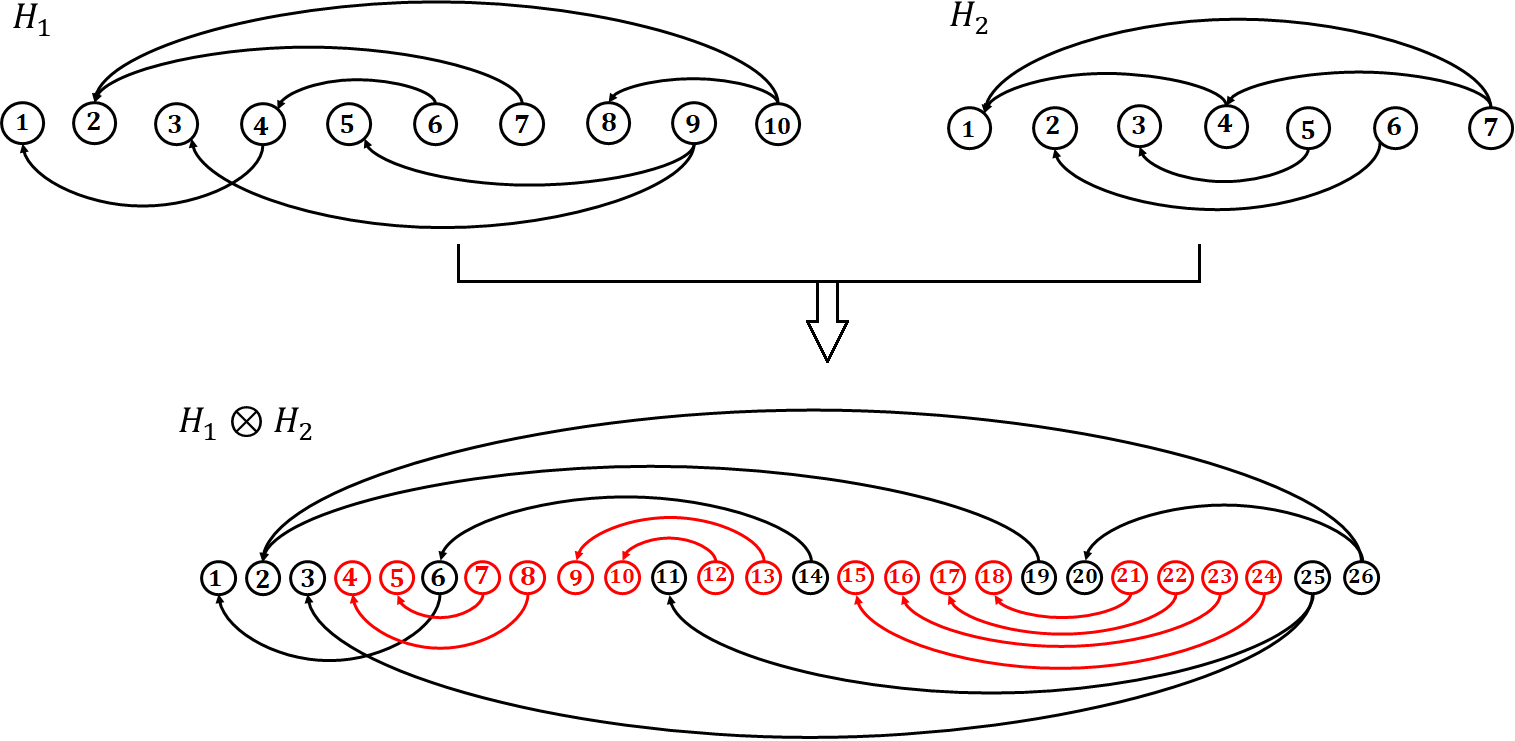}
	\caption{Super nebula $H_{1}$ drawn under its super nebula ordering $\theta_{1}$, $\Delta$galaxy $H_{2}$ drawn under its $\Delta$galaxy ordering $\theta_{2}$, and a corresponding key tournament $H_{1}\otimes H_{2}$ under $\theta$. All the non drawn arcs are forward.}
	\label{fig:keytour}
\end{figure}

\begin{figure}[h]
	\centering
	\includegraphics[width=1\linewidth]{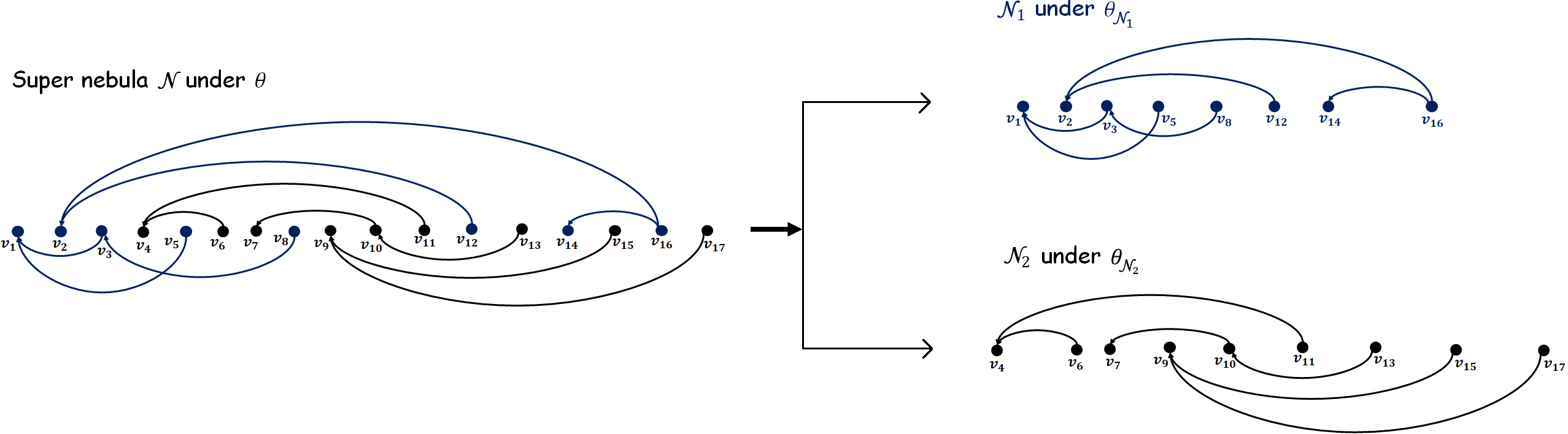}
	\caption{Super nebula $\mathcal{N}$ under $\theta$, $\mathcal{N}_{1}$ under $\theta_{1}$, and $\mathcal{N}_{2}$ under $\theta_{2}$. All the arcs that are not drawn are forward.}
	\label{fig:separate}
\end{figure}

Let $\mathcal{G}$ be a regular $\Delta$galaxy under $\alpha=(z_{1},...,z_{n_{2}})$ and let $\Delta =\lbrace z_{j_{1}},z_{j_{2}},z_{j_{3}}\rbrace$ be the triangle of $\mathcal{G}$ under $\alpha$. Define $G_{\mathcal{G},\alpha}=\lbrace \tilde{\mathcal{G}}; \tilde{\mathcal{G}}$ is obtained from $\mathcal{G}$ by reversing the orientation of exactly one arc of $\Delta \rbrace$. Notice that $\mid$$G_{\mathcal{G},\alpha}$$\mid$ $=3$. Let $\mathcal{G}_{\overline{\Delta}}= \mathcal{G}\backslash \lbrace z_{j_{1}},z_{j_{2}},z_{j_{3}}\rbrace$ and let $\alpha_{\overline{\Delta}}$ be the restriction of $\alpha$ to $V(\mathcal{G}_{\overline{\Delta}})$.\\
Let $\mathcal{N}$ be a regular super nebula under $\theta = (v_{1},...,v_{n_{1}})$ and let $B^{\mathcal{N},\theta}=\lbrace t_{b}^{i}=(v_{i_{1}}, v_{i_{2}},v_{i_{3}}); i=1,...,s\rbrace$ (note that for all $i\in\lbrace 1,...,s\rbrace$, $i_{1}<i_{2}<i_{3}$). Notice that $\forall 1\leq i\leq s$, $ \mid$$E(B(\mathcal{N}$$\mid$$\lbrace v_{i_{1}},v_{i_{2}},v_{i_{3}}\rbrace$$,(v_{i_{1}},v_{i_{2}},v_{i_{3}})))$$\mid =2$. $\forall 1\leq i\leq s$, let $e_{i}\in E(\mathcal{N}$$\mid$$\lbrace v_{i_{1}},v_{i_{2}},v_{i_{3}}\rbrace)$ such that $e_{i}$ is forward under $\theta$. $e_{i}$ is called the \textit{forward arc of $(v_{i_{1}},v_{i_{2}},v_{i_{3}})$}. Define $E^{\mathcal{N},\theta}_{f}=\lbrace e_{1},...,e_{s}\rbrace$. The \textit{mutant super nebula} $\tilde{\mathcal{N}}$ under $\theta$ is the digraph obtained from $\mathcal{N}$ by deleting all the arcs in $E^{\mathcal{N},\theta}_{f}$.\\
A regular super nebula $\mathcal{K}=\mathcal{N}\otimes\mathcal{G}$ under $\tilde{\theta}$ is a \textit{key tournament corresponding to $\mathcal{N}$ and $\mathcal{G}$ under $\theta$ and $\alpha$ respectively} if $\mathcal{K}$ under $\tilde{\theta}$ satisfies all the following:
\begin{itemize}
\item $V(\mathcal{K})=V(\mathcal{N})\cup \bigcup_{i=1}^{s}U_{i}$ where $U_{i}=\lbrace u_{1}^{i},...,u_{n_{2}-3}^{i}\rbrace$ and $s=\mid$$B^{\mathcal{N},\theta}$$\mid$.
\item $\mathcal{K}$$\mid$$V(\mathcal{N})$ is isomorphic to $\mathcal{N}$ and the restriction of $\tilde{\theta}$ to $V(\mathcal{N})$ is the super nebula ordering $\theta$ of $\mathcal{N}$. 
\item $\mathcal{K}$$\mid$$U_{i}$ is the galaxy tournament $\mathcal{G}_{\overline{\Delta}}$ and the restriction of $\tilde{\theta}$ to $U_{i}$  (say $\tilde{\theta}_{U_{i}}$) is the galaxy ordering $\alpha_{\overline{\Delta}}$ of $\mathcal{G}_{\overline{\Delta}}$ for $i=1,...,s$.
\item $B^{\mathcal{K},\tilde{\theta}}=B^{\mathcal{N},\theta}=\lbrace (v_{i_{1}}, v_{i_{2}},v_{i_{3}}); i=1,...,s\rbrace$.
\item $\mathcal{K}$$\mid$$X_{i}$ is isomorphic to $\tilde{\mathcal{G}}$ for some $\tilde{\mathcal{G}}\in G_{\mathcal{G},\alpha}$, where $X_{i}=\lbrace v_{i_{1}}, v_{i_{2}},v_{i_{3}}\rbrace\cup U_{i}$ for $i=1,...,s$. Moreover, $\tilde{\theta}_{X_{i}}$,  the restriction of $\tilde{\theta}$ to $X_{i}$ for $i=1,...,s$ verifies the following: let $e_{i}\in \mathcal{K}$$\mid$$\lbrace v_{i_{1}}, v_{i_{2}},v_{i_{3}}\rbrace$ such that $e_{i}$ is forward under $\tilde{\theta}$. The tournament obtained from $\mathcal{K}$$\mid$$X_{i}$ by reversing the orientation of $e_{i}$ is the tournament $\mathcal{G}$ and the ordering $\tilde{\theta}_{X_{i}}$ is then the triangular ordering $\alpha$ of $\mathcal{G}$. Here we say that $U_{i}$ \textit{corresponds to} $t_{i} \in B^{\mathcal{K},\tilde{\theta}}$ for $i=1,...,s$. 
\item The backward arcs of $\mathcal{K}$ under $\tilde{\theta}$ are exactly the backward arcs of $\mathcal{K}$$\mid$$V(\mathcal{N})$ under $\theta$ and the backward arcs of $\mathcal{K}$$\mid$$U_{i}$ under $\tilde{\theta}_{U_{i}}$ for $i=1,...,s$.
 \end{itemize}
 In Figure \ref{fig:keytour} we draw a super nebula $H_{1}$ under $\theta_{1}=(1,...,10)$, a $\Delta$galaxy $H_{2}$ under $\theta_{2}=(\mathbf{1,...,7})$, and a key tournament $H_{1}\otimes H_{2}$ under $\theta =(\mathit{1,...,26})$ corresponding to $H_{1}$ and $ H_{2}$ under $ \theta_{1} $ and $ \theta_{2} $ respectively (note that we have more than one key tournament corresponding to $H_{1}$ and $ H_{2}$ under $ \theta_{1} $ and $ \theta_{2} $ respectively).\vspace{3mm}
\\Let $\mathcal{N}$ be a regular super nebula under $\theta =(v_{1},...,v_{n})$. Let $Q_{1},...,Q_{m}$ be the stars of $\mathcal{N}$ under $\theta$ and let $\Sigma_{1},...,\Sigma_{l}$ be the super $2$-nebulas of $\mathcal{N}$ under $\theta$. Let $\mathcal{N}_{1}=\mathcal{N}$$\mid$$\bigcup_{i=1}^{l}V(\Sigma_{i})$ and let $\mathcal{N}_{2}=\mathcal{N}$$\mid$$\bigcup_{i=1}^{m}V(Q_{i})$.  For $k_{1}\in \lbrace 0,...,l\rbrace$ define  $\mathcal{N}_{1}^{k_{1}}= \mathcal{N}$$\mid$$\bigcup_{i=1}^{k_{1}}V(\Sigma_{i})$ and let $\theta_{1}^{k_{1}}$ be the restriction of $\theta$ to $V(\mathcal{N}_{1}^{k_{1}})$, where $\mathcal{N}_{1}^{l} = \mathcal{N}_{1}$, and $\mathcal{N}_{1}^{0}$ is the empty tournament. For $k_{2}\in \lbrace 0,...,m\rbrace$ define $\mathcal{N}^{k_{2}}_{2}= \mathcal{N}$$\mid$$\bigcup_{i=1}^{k_{2}}V(Q_{i})$ and let $\theta_{2}^{k_{2}}$ be the restriction of $\theta$ to $V(\mathcal{N}_{2}^{k_{2}})$, where $\mathcal{N}_{2}^{m} = \mathcal{N}_{2}$, and $\mathcal{N}_{2}^{0}$ is the empty tournament. For $i\in \lbrace 1,2\rbrace$,  let $s^{\mathcal{N},\theta}_{\mathcal{N}^{k_{i}}_{i}}$ be the restriction of $s^{\mathcal{N},\theta}$ to the $0's$ and $1's$ corresponding to $V(\mathcal{N}^{k_{i}}_{i})$ (notice that $s^{\mathcal{N},\theta}_{\mathcal{N}^{k_{i}}_{i}}= s^{\mathcal{N}^{k_{i}}_{i},\theta_{i}^{k_{i}}}$) and let $^{c}s^{\mathcal{N},\theta}_{\mathcal{N}^{k_{i}}_{i}}$ be the vector obtained from $s^{\mathcal{N},\theta}_{\mathcal{N}^{k_{i}}_{i}}$ by replacing every subsequence of consecutive $1's$ corresponding to the same entry of $s^{\mathcal{N},\theta}_{c}$ by single $1$.\\ In the example in Figure \ref{fig:separate} we have: $s^{\mathcal{N},\theta}=(0,0,0,0,1,1,1,1,0,0,1,1,1,1,1,0,1)$, $s^{\mathcal{N},\theta}_{c}=(0,0,0,0,1,0,0,1,0,1)$, $s^{\mathcal{N},\theta}_{\mathcal{N}_{1}}=(0,0,0,1,1,1,1,0)$, $s^{\mathcal{N},\theta}_{\mathcal{N}_{2}}=(0,1,1,0,0,1,1,1,1)$, $^{c}s^{\mathcal{N},\theta}_{\mathcal{N}_{1}}=(0,0,0,1,1,0)$, $^{c}s^{\mathcal{N},\theta}_{\mathcal{N}_{2}}=(0,1,0,0,1,1)$.
\\We say that a smooth $(c,\lambda ,w)$$-$structure of a tournament $T$ \textit{corresponds to $\mathcal{N}$ under $\theta$} if $w = s_{c}^{\mathcal{N},\theta}$. For $i\in \lbrace 1,2\rbrace$, we say that a smooth $(c,\lambda ,w)$$-$structure of a tournament $T$ \textit{corresponds to $\mathcal{N}^{k_{i}}_{i}$ under $(\mathcal{N},\theta )$} if $w= ^{c}s^{\mathcal{N},\theta}_{\mathcal{N}^{k_{i}}_{i}}$.\\
Let $\delta^{s_{c}^{\mathcal{N},\theta}}:$ $\lbrace j: s_{c_{j}}^{\mathcal{N},\theta} = 1 \rbrace \rightarrow \mathbb{N}$ be a function that assigns to every nonzero entry of $s_{c}^{\mathcal{N},\theta}$ the number of consecutive $1'$s of $s^{\mathcal{N},\theta}$ replaced by that entry of $s_{c}^{\mathcal{N},\theta}$. Similarly let $\delta^{\nu^{i}}:$ $\lbrace j: \nu^{i}_{j} = 1 \rbrace \rightarrow \mathbb{N}$ be a function that assigns to every nonzero entry of $\nu^{i}$ the number of consecutive $1'$s of $s^{\mathcal{N},\theta}_{\mathcal{N}_{i}^{k_{i}}}$ replaced by that entry of $\nu^{i}$ for $i=1,2$, where $\nu^{i}=^{c}s^{\mathcal{N},\theta}_{\mathcal{N}^{k_{i}}_{i}}$ for $i=1,2$ (recall that: $k_{1}\in \lbrace 0,...,l\rbrace$, and $k_{2}\in \lbrace 0,...,m\rbrace$). \vspace{2mm}\\ 
Let $\mathcal{N}$ be a regular super nebula under $\theta =(v_{1},...,v_{n})$. Let $Q_{1},...,Q_{m}$ be the stars of $\mathcal{N}$ under $\theta$ and let $\Sigma_{1},...,\Sigma_{l}$ be the super $2$-nebulas of $\mathcal{N}$ under $\theta$. Let $\mathcal{N}_{1}=\mathcal{N}$$\mid$$\bigcup_{i=1}^{l}V(\Sigma_{i})$ and let $\mathcal{N}_{2}=\mathcal{N}$$\mid$$\bigcup_{i=1}^{m}V(Q_{i})$. Fix $k\in\lbrace 0,...,l\rbrace$ and $z\in\lbrace 0,...,m\rbrace$. Let $\tilde{\mathcal{N}}^{k}_{1}= \tilde{\mathcal{N}}$$\mid$$V(\mathcal{N}^{k}_{1})$ and let $\tilde{\mathcal{N}}^{z}_{2}= \tilde{\mathcal{N}}$$\mid$$V(\mathcal{N}^{z}_{2})$, where $\tilde{\mathcal{N}}$ is the mutant super nubula obtained from $\mathcal{N}$ under $\theta$. Let $\theta_{1}^{k}=(v_{k_{1}},...,v_{k_{q_{k}}})$ and $\theta_{2}^{z}=(v_{z_{1}},...,v_{z_{p_{z}}})$. Let $(S_{1},...,S_{\mid w \mid})$ be a smooth $(c,\lambda ,w)$$-$structure corresponding to $\mathcal{N}$ under $\theta$ (resp. corresponding to $\mathcal{N}^{k}_{1}$ under $(\mathcal{N},\theta )$) (resp. corresponding to $\mathcal{N}^{z}_{2}$ under $(\mathcal{N},\theta )$). Let $i_{r}$ be such that $w(i_{r}) = 1$. Assume that $S_{i_{r}} = \lbrace s^{1}_{i_{r}},...,s_{i_{r}}^{\mid S_{i_{r}} \mid} \rbrace$ and $(s^{1}_{i_{r}},...,s_{i_{r}}^{\mid S_{i_{r}} \mid})$ is a transitive ordering. Write $m(i_{r}) = \lfloor\frac{\mid S_{i_{r}} \mid}{\delta^{w}(i_{r})}\rfloor$.\\ Denote $S^{j}_{i_{r}} = \lbrace s^{(j-1)m(i_{r})+1}_{i_{r}},...,s_{i_{r}}^{jm(i_{r})} \rbrace$ for $j \in \lbrace 1,...,\delta^{w}(i_{r}) \rbrace$. For every $v \in S^{j}_{i_{r}}$ denote $\xi(v) = (\mid$$\lbrace k < i_{r}: w(k) = 0 \rbrace$$\mid$ $+$ $\displaystyle{\sum_{k < i_{r}: w(k) = 1}\delta^{w}(k) })$ $+$ $j$. For every $v \in S_{i_{r}}$ such that $w(i_{r}) = 0$ denote $\xi(v) = (\mid$$\lbrace k < i_{r}: w(k) = 0 \rbrace$$\mid$ $+$ $\displaystyle{\sum_{k < i_{r}: w(k) = 1}\delta^{w}(k) })$ $+$ $1$. We say that $\tilde{\mathcal{N}}$ (resp. $\tilde{\mathcal{N}_{1}^{k}}$) (resp. $\tilde{\mathcal{N}^{z}_{2}}$) is \textit{well-contained in} $(S_{1},...,S_{\mid w \mid})$ that corresponds  to $\mathcal{N}$ under $\theta$ (resp. corresponds to $\mathcal{N}^{k}_{1}$ under $(\mathcal{N},\theta )$) (resp. corresponds to $\mathcal{N}^{z}_{2}$ under $(\mathcal{N},\theta )$) if there is an injective homomorphism $f$ of $\tilde{\mathcal{N}}$ (resp. $\tilde{\mathcal{N}^{k}_{1}}$) (resp. $\tilde{\mathcal{N}^{z}_{2}}$) into $T$$\mid$$\bigcup_{i = 1}^{\mid w \mid}S_{i}$ such that $\xi(f(v_{j})) = j$ for every $j \in \lbrace 1,...,n \rbrace$ (resp. $\xi(f(v_{k_{j}})) = j$ for every $j \in \lbrace 1,...,q_{k} \rbrace$) (resp. $\xi(f(v_{z_{j}})) = j$ for every $j \in \lbrace 1,...,p_{z} \rbrace$). 
\subsection{Proof of Theorem \ref{p}}
We start by the following technical lemma:
\begin{lemma}\label{m}
Let $\mathcal{N}$ be a regular super nebula under $\theta_{1}$ with $\mid$$\mathcal{N}$$\mid =\mu_{1}$ and let $\mathcal{G}$ be a regular $\Delta$galaxy under $\theta_{2}$. Let $\mathcal{K}=\mathcal{N}\otimes\mathcal{G}$ under $\theta$ be a key tournament corresponding to $\mathcal{N}$ and $\mathcal{G}$ under $\theta_{1}$ and $\theta_{2}$ respectively. Let $Q_{1},...,Q_{m}$ be the stars of $\mathcal{K}$ under $\theta$ and let $\Sigma_{1},...,\Sigma_{l}$ be the super $2$-nebulas of $\mathcal{K}$ under $\theta$. Let $\mathcal{K}_{1}=\mathcal{K}$$\mid$$\bigcup_{i=1}^{l}V(\Sigma_{i})$ and let $\mathcal{K}_{2}=\mathcal{K}$$\mid$$\bigcup_{i=1}^{m}V(Q_{i})$. Let $0 < \lambda < \frac{1}{(2\mu_{1})^{\mu_{1}+2}}$, $c > 0$ be constants, and $w$ be a $\lbrace 0,1 \rbrace$$-$vector. Fix $k \in \lbrace 0,...,l \rbrace$ and let $\widehat{\lambda} = (2\mu_{1})^{l-k}\lambda$ and $\widehat{c} = \frac{c}{(2\mu_{1})^{l-k}}$. There exist $ \epsilon_{k} > 0$ such that $\forall 0 < \epsilon < \epsilon_{k}$, for every $\epsilon$$-$critical tournament $T$ with $\mid$$T$$\mid$ $= n$ containing $\chi = (S_{1},...,S_{\mid w \mid})$ as a smooth $(\widehat{c},\widehat{\lambda},w)$$-$structure corresponding to $\mathcal{K}_{1}^{k}$ under $(\mathcal{K},\theta )$, we have $\tilde{\mathcal{K}}_{1}^{k}$ is well-contained in $\chi$.  
\end{lemma}
\begin{proof}
The proof is by induction on $k$. For $k=0$ the statement is obvious since $\tilde{\mathcal{K}}^{0}_{1}$ is the empty digraph. Suppose that $\chi = (S_{1},...,S_{\mid w \mid})$ is a smooth $(\widehat{c},\widehat{\lambda},w)$$-$structure in $T$ corresponding to $\mathcal{K}_{1}^{k}$ under $(\mathcal{K},\theta )$ with $\theta = (h_{1},...,h_{\mid \mathcal{K} \mid})$ and $\mid $$\mathcal{K}$$ \mid = h$. Let $\theta_{1}^{k}=(h_{p_{1}},...,h_{p_{q}})$ be the restriction of $\theta$ to $V(\mathcal{K}_{1}^{k})$ (notice that $q\leq \mu_{1}$). Let $h_{p_{b_{0}}}$ and $h_{p_{b_{1}}}$  be the center of $\Sigma_{k}$, and let $h_{p_{b_{2}}},...,h_{p_{b_{r}}},...,h_{p_{b_{z}}}$ be its leaves for some integer $z>0$ such that $h_{p_{b_{2}}},...,h_{p_{b_{r}}}$ are the leaves incident to $h_{p_{b_{0}}}$ and $h_{p_{b_{r+1}}},...,h_{p_{b_{z}}}$ are the leaves incident to $h_{p_{b_{1}}}$ (note that $\forall 2\leq i<j\leq z$, we don't have necessarily $b_{i}<b_{j}$). 
$\forall 0 \leq i \leq z$, let $D_{i} = \lbrace v \in \displaystyle{\bigcup_{j=1}^{\mid w \mid}S_{j}};$ $ \xi(v) = b_{i} \rbrace$.
Then $\exists x \in \lbrace 1,...,\mid$$w$$\mid \rbrace$, $\exists y \in \lbrace 1,...,\mid$$w$$\mid \rbrace$, $\exists f \in \lbrace 1,...,\mid$$w$$\mid \rbrace$ with $x<y$, $w(x)=w(y)=0$, and $w(f)=1$, such that $D_{0} = S_{x}$, $D_{1} = S_{y}$, and $\forall 2 \leq i \leq z$, $D_{i} \subseteq S_{f}$. Since we can assume that $\epsilon < min \lbrace log_{\frac{\widehat{c}}{2\mu_{1}}}(1-\frac{\widehat{c}}{\mu_{1}}),log_{\frac{\widehat{c}}{4}}(\frac{1}{2})\rbrace$, then by Lemma \ref{s} there exists vertices $d_{0},d_{1},...,d_{z}$ such that $d_{i} \in D_{i}$ for $i=0,1,...,z$ and \\
$\ast$ $d_{0}\leftarrow \lbrace d_{1},d_{2},...,d_{r}\rbrace$ and $\lbrace d_{r+1},...,d_{z}\rbrace \leftarrow d_{1}$ if $x<f<y$.\\
$\ast$ $\lbrace d_{0},d_{r+1},...,d_{z}\rbrace \leftarrow d_{1}$ and $\lbrace d_{2},...,d_{r}\rbrace \leftarrow d_{0}$  if $f<x<y$.\\
$\ast$ $d_{0}\leftarrow \lbrace d_{1},d_{2},...,d_{r}\rbrace$ and $d_{1}\leftarrow \lbrace d_{r+1},...,d_{z}\rbrace$  if $x<y<f$.\\
So $T$$\mid$$\lbrace d_{0},d_{1},d_{2},...,d_{z} \rbrace$ contains a copy of $\tilde{\mathcal{K}_{1}^{k}}$$\mid$$V(\Sigma_{k})$. Denote this copy by $Y$.\\    
$\forall i \in \lbrace 1,...,\mid$$w$$\mid \rbrace \backslash \lbrace x,y,f \rbrace$, let $S_{i}^{*} = \displaystyle{\bigcap_{p\in V(Y)}S_{i,p}}$.  Then by Lemma \ref{g}, $\mid$$S_{i}^{*}$$\mid$ $\geq (1-\mid$$Y$$\mid\widehat{\lambda})\mid$$S_{i}$$\mid$ $\geq (1-\mu_{1}\widehat{\lambda})\mid$$S_{i}$$\mid$ $\geq \frac{1}{2\mu_{1}}\mid$$S_{i}$$\mid$ since $\widehat{\lambda} \leq \frac{2\mu_{1}-1}{2\mu_{1}^{2}}$.
Write $\mathcal{H} = \lbrace 1,...,q \rbrace \backslash \lbrace b_{0},...,b_{z} \rbrace$. If $\lbrace v\in S_{f}: \xi(v) \in \mathcal{H} \rbrace \neq \phi$, then define $J_{f} = \lbrace \eta \in \mathcal{H}: \exists v \in S_{f}$ and $\xi(v)= \eta \rbrace$. Now $\forall \eta \in J_{f}$, let $S_{f}^{*\eta}= \lbrace v \in S_{f}: \xi(v)=\eta$ and $v \in \displaystyle{\bigcap_{i\in \lbrace 0,1\rbrace}S_{f,d_{i}}} \rbrace$. Then by Lemma \ref{g}, $\forall \eta \in J_{f}$, we have $\mid$$S_{f}^{*\eta}$$\mid$ $ \geq \frac{1-2\mu_{1}\widehat{\lambda}}{\mu_{1}}\mid $$S_{f}$$\mid $ $\geq \frac{\mid S_{f}\mid}{2\mu_{1}}$ since $\widehat{\lambda} \leq \frac{1}{4\mu_{1}}$. Now $\forall \eta \in J_{f}$, select arbitrary $\lceil \frac{\mid S_{f}\mid}{2\mu_{1}}\rceil$ vertices of $S_{f}^{*\eta}$ and denote the union of these $\mid$$J_{f}$$\mid$ sets by $S_{f}^{*}$. 
So we have defined $t$ sets $S_{1}^{*},...,S^{*}_{t}$, where $t = \mid$$w$$\mid$ $-2$ if $S_{f}^{*}$ is defined and $t = \mid$$w$$\mid$ $-3$ if $S_{f}^{*}$ is not defined. We have $\mid$$S_{i}^{*}$$\mid$ $\geq \frac{\widehat{c}}{2\mu_{1}}tr(T)$ for every defined $S_{i}^{*}$ with $w(i) = 1$, and $\mid$$S_{i}^{*}$$\mid$ $\geq \frac{\widehat{c}}{2\mu_{1}}n$ for every defined $S_{i}^{*}$ with $w(i) = 0$. Now Lemma \ref{b} implies that $\chi^{*}=(S_{1}^{*},...,S^{*}_{t})$ form a smooth $(\frac{\widehat{c}}{2\mu_{1}},2\mu_{1}\widehat{\lambda},w^{*})$$-$structure of $T$ corresponding to $\mathcal{K}_{1}^{k-1}$ under $(\mathcal{K},\theta )$, where $\frac{\widehat{c}}{2\mu_{1}}= \frac{c}{(2\mu_{1})^{l-(k-1)}}, 2\mu_{1}\widehat{\lambda}=(2\mu_{1})^{l-(k-1)}\lambda$, and $w^{*}$ is an appropriate $\lbrace 0,1 \rbrace$$-$vector.
Now take $\epsilon_{k} < min \lbrace \epsilon_{k-1}, log_{\frac{\widehat{c}}{2\mu_{1}}}(1-\frac{\widehat{c}}{\mu_{1}}),log_{\frac{\widehat{c}}{4}}(\frac{1}{2}) \rbrace$. So by induction hypothesis $\tilde{\mathcal{K}}_{1}^{k-1}$ is well-contained in $\chi^{*}$. Now by merging the well-contained copy of $\tilde{\mathcal{K}}_{1}^{k-1}$ and $Y$ we get a well-contained copy of $\tilde{\mathcal{K}}_{1}^{k}$. $\blacksquare$ 
\end{proof}
\vspace{3mm}\\We also need the following technical lemma:
\begin{lemma}\label{o}
Let $\mathcal{N}$ be a regular super nebula under $\theta_{1}$ with $\mid$$\mathcal{N}$$\mid =\mu_{1}$ and let $\mathcal{G}$ be a regular $\Delta$galaxy under $\theta_{2}$ with $\mid$$\mathcal{G}$$\mid =\mu_{2}$. Let $\delta =\mu_{1}\mu_{2}$. Let $\mathcal{K}=\mathcal{N}\otimes\mathcal{G}$ under $\theta$ be a key tournament corresponding to $\mathcal{N}$ and $\mathcal{G}$ under $\theta_{1}$ and $\theta_{2}$ respectively. Let $Q_{1},...,Q_{m}$ be the stars of $K$ under $\theta$ and let $\Sigma_{1},...,\Sigma_{l}$ be the super $2$-nebulas of $K$ under $\theta$. Let $\mathcal{K}_{1}=\mathcal{K}$$\mid$$\bigcup_{i=1}^{l}V(\Sigma_{i})$ and let $\mathcal{K}_{2}=\mathcal{K}$$\mid$$\bigcup_{i=1}^{m}V(Q_{i})$. Let $0 < \lambda < \frac{1}{(2\delta)^{\delta+2}}$, $c > 0$ be constants, and $w$ be a $\lbrace 0,1 \rbrace$$-$vector. Fix $k \in \lbrace 0,...,m \rbrace$ and let $\widehat{\lambda} = (2\delta)^{m-k}\lambda$ and $\widehat{c} = \frac{c}{(2\delta)^{m-k}}$. There exist $ \epsilon_{k} > 0$ such that $\forall 0 < \epsilon < \epsilon_{k}$, for every $\epsilon$$-$critical tournament $T$ with $\mid$$T$$\mid$ $= n$ containing $\chi = (S_{1},...,S_{\mid w \mid})$ as a smooth $(\widehat{c},\widehat{\lambda},w)$$-$structure corresponding to $\mathcal{K}_{2}^{k}$ under $(\mathcal{K},\theta )$, we have $\tilde{\mathcal{K}}_{2}^{k}$ is well-contained in $\chi$.
\end{lemma}
\begin{proof}
The proof is by induction on $k$. For $k=0$ the statement is obvious since $\tilde{\mathcal{K}}_{2}^{0}$ is the empty digraph. Suppose that $\chi = (S_{1},...,S_{\mid w \mid})$ is a smooth $(\widehat{c},\widehat{\lambda},w)$$-$structure in $T$ corresponding to $\mathcal{K}_{2}^{k}$ under $(\mathcal{K},\theta )$ with $\theta = (h_{1},...,h_{\mid \mathcal{K} \mid})$ and $\mid $$\mathcal{K}$$ \mid = h$. Let $\theta_{2}^{k}=(h_{q_{1}},...,h_{q_{p}})$ be the restriction of $\theta$ to $V(\mathcal{K}_{2}^{k})$ (notice that $p\leq \delta$). Let $h_{q_{a_{0}}}$  be the center of $Q_{k}$, and let $h_{q_{a_{1}}},...,h_{q_{a_{d}}}$ be its leaves for some integer $d>0$.
$\forall 0 \leq i \leq d$, let $R_{i} = \lbrace v \in \displaystyle{\bigcup_{j=1}^{\mid w \mid}S_{j}};$ $ \xi(v) = a_{i} \rbrace$. We have $\mid$$R_{0}$$\mid$ $\geq \widehat{c}n$ and $R_{i}\geq\frac{\widehat{c}}{\delta}tr(T)$ for $i=1,...,d$. Since we can assume that $\epsilon < log_{\frac{\widehat{c}}{2\delta}}(1-\frac{\widehat{c}}{\delta})$, then by Lemma \ref{r} there exists vertices $r_{0},r_{1},...,r_{d}$ such that $r_{i} \in R_{i}$ for $i=0,1,...,d$ and \vspace{3mm}\\
$\ast$ $r_{0}\leftarrow \lbrace r_{1},...,r_{d}\rbrace$ if $Q_{k}$ is a left star of $\mathcal{K}$ under $\theta$.\\
$\ast$ $\lbrace r_{1},...,r_{d}\rbrace \leftarrow r_{0}$ if $Q_{k}$ is a right star of $\mathcal{K}$ under $\theta$.\\
$\ast$ $\lbrace r_{1},...,r_{d_{1}}\rbrace \leftarrow r_{0}\leftarrow \lbrace r_{d_{1}+1},...,r_{d}\rbrace$ if $Q_{k}$ is a middle star of $\mathcal{K}$ under $\theta$.\vspace{3mm}\\
So $T$$\mid$$\lbrace r_{0},r_{1},...,r_{d} \rbrace$ contains a copy of $\tilde{\mathcal{K}}^{k}_{2}$$\mid$$V(Q_{k})$. Denote this copy by $Y$. Let $ x \in \lbrace 1,...,\mid$$w$$\mid \rbrace$ such that $R_{0}= S_{x}$ and let $y_{1},...,y_{d} \in \lbrace 1,...,\mid$$w$$\mid \rbrace\backslash\lbrace x\rbrace$ such that $R_{i}\subseteq S_{y_{i}}$ for $i=1,...,d$. Notice that we have: $w(x)=0$ and $w(y_{i})=1$ for $i=1,...,d$. Also notice that we don't have necessarily that $y_{1},...,y_{d}$ are distinct. 
    \\$\forall i \in \lbrace 1,...,\mid$$w$$\mid \rbrace \backslash \lbrace x,y_{1},...,y_{d} \rbrace$, let $S_{i}^{*} = \displaystyle{\bigcap_{j=0}^{d}S_{i,r_{j}}}$.  Then by Lemma \ref{g}, $\mid$$S_{i}^{*}$$\mid$ $\geq (1-\mid$$Y$$\mid\widehat{\lambda})\mid$$S_{i}$$\mid$ $\geq (1-\delta\widehat{\lambda})\mid$$S_{i}$$\mid$ $\geq \frac{1}{2\delta}\mid$$S_{i}$$\mid$ since $\widehat{\lambda} \leq \frac{2\delta-1}{2\delta^{2}}$.
Write $\mathcal{H} = \lbrace 1,...,p \rbrace \backslash \lbrace a_{0},...,a_{d} \rbrace$. Let $Y_{i}=\lbrace v\in V(Y): v\in S_{y_{i}}\rbrace$ for $i=1,...,d$. $\forall 1\leq i\leq d$, if $\lbrace v\in S_{y_{i}}: \xi(v) \in \mathcal{H} \rbrace \neq \phi$, then define $J_{y_{i}} = \lbrace \eta \in \mathcal{H}: \exists v \in S_{y_{i}}$ and $\xi(v)= \eta \rbrace$. Now $\forall \eta \in J_{y_{i}}$, let $S_{y_{i}}^{*\eta}= \lbrace v \in S_{y_{i}}: \xi(v)=\eta$ and $v \in \displaystyle{\bigcap_{q\in V(Y)\backslash Y_{i}}S_{y_{i},q}} \rbrace$. Then by Lemma \ref{g}, $\forall \eta \in J_{y_{i}}$, we have $\mid$$S_{y_{i}}^{*\eta}$$\mid$ $\geq (1-(d-1)\delta\widehat{\lambda})\frac{\mid S_{y_{i}}\mid}{\delta} \geq \frac{1-\delta^{2}\widehat{\lambda}}{\delta}\mid $$S_{y_{i}}$$\mid $ $\geq \frac{\mid S_{y_{i}}\mid}{2\delta}$ since $\widehat{\lambda} \leq \frac{1}{2\delta^{2}}$. Now $\forall \eta \in J_{y_{i}}$, select arbitrary $\lceil \frac{\mid S_{y_{i}}\mid}{2\delta}\rceil$ vertices of $S_{y_{i}}^{*\eta}$ and denote the union of these $\mid$$J_{y_{i}}$$\mid$ sets by $S_{y_{i}}^{*}$. 
So we have defined some number of sets. Denote by $t$ the number of these defined sets and by $S_{1}^{*},...,S^{*}_{t}$ these  sets. We have $\mid$$S_{i}^{*}$$\mid$ $\geq \frac{\widehat{c}}{2\delta}tr(T)$ for every defined $S_{i}^{*}$ with $w(i) = 1$, and $\mid$$S_{i}^{*}$$\mid$ $\geq \frac{\widehat{c}}{2\delta}n$ for every defined $S_{i}^{*}$ with $w(i) = 0$. Now Lemma \ref{b} implies that $\chi^{*}=(S_{1}^{*},...,S^{*}_{t})$ form a smooth $(\frac{\widehat{c}}{2\delta},2\delta\widehat{\lambda},w^{*})$$-$structure of $T$ corresponding to $\mathcal{K}_{2}^{k-1}$ under $(\mathcal{K},\theta )$, where $\frac{\widehat{c}}{2\delta}= \frac{c}{(2\delta)^{l-(k-1)}}, 2\delta\widehat{\lambda}=(2\delta)^{l-(k-1)}\lambda$, and $w^{*}$ is an appropriate $\lbrace 0,1 \rbrace$$-$vector.
Now take $\epsilon_{k} < min \lbrace \epsilon_{k-1}, log_{\frac{\widehat{c}}{2\delta}}(1-\frac{\widehat{c}}{\delta}) \rbrace$. So by induction hypothesis $\tilde{\mathcal{K}}^{k-1}_{2}$ is well-contained in $\chi^{*}$. Now by merging the well-contained copy of $\tilde{\mathcal{K}}^{k-1}_{2}$ and $Y$ we get a well-contained copy of $\tilde{\mathcal{K}^{k}_{2}}$. $\blacksquare$\vspace{4mm} 
\end{proof}

From Lemmas \ref{m} and \ref{o} we get the following lemma:
\begin{lemma}
Let $\mathcal{N}$ be a regular super nebula under $\theta_{1}$ with $\mid$$\mathcal{N}$$\mid =\mu_{1}$ and let $\mathcal{G}$ be a regular $\Delta$galaxy under $\theta_{2}$ with $\mid$$\mathcal{G}$$\mid =\mu_{2}$. Let $\delta =\mu_{1}\mu_{2}$. Let $\mathcal{K}=\mathcal{N}\otimes\mathcal{G}$ under $\theta =(h_{1},...,h_{h})$ be a key tournament corresponding to $\mathcal{N}$ and $\mathcal{G}$ under $\theta_{1}$ and $\theta_{2}$ respectively.  Let $0 < \lambda_{0} < \frac{1}{(2\delta)^{\delta+3}}$, $c > 0$ be constants, and $w$ be a $\lbrace 0,1 \rbrace$$-$vector. There exist $ \epsilon_{0} > 0$ such that for every $\epsilon$$-$critical tournament $T$ with $\mid$$T$$\mid$ $= n$ containing $\chi = (A_{1},...,A_{\mid w \mid})$ as a smooth $(c_{0},\lambda_{0},w)$$-$structure corresponding to $\mathcal{K}$ under the ordering $\theta$, we have $\tilde{\mathcal{K}}$ is well-contained in $\chi$.  
\end{lemma}
\begin{proof}
Let $Q_{1},...,Q_{m}$ be the stars of $\mathcal{K}$ under $\theta$ and let $\Sigma_{1},...,\Sigma_{l}$ be the super $2$-nebulas of $\mathcal{K}$ under $\theta$. Let $\mathcal{K}_{1}=\mathcal{K}$$\mid$$\bigcup_{i=1}^{l}V(\Sigma_{i})$ and let $\mathcal{K}_{2}=\mathcal{K}$$\mid$$\bigcup_{i=1}^{m}V(Q_{i})$. Let $\theta_{1} =(h_{p_{1}},...,h_{p_{q}})$ be the restriction of $\theta$ to $V(\mathcal{K}_{1})$ and let $\theta_{2} =(h_{q_{1}},...,h_{q_{p}})$ be the restriction of $\theta$ to $V(\mathcal{K}_{2})$. $\forall 1\leq i \leq q$, let  $S_{i} = \lbrace v \in \bigcup_{j=1}^{\mid w \mid}A_{j};$ $ \xi(v) = p_{i} \rbrace$. $\forall 1\leq j \leq $ $\mid$$w$$\mid$, let $S^{*}_{j}=  \displaystyle{\bigcup_{S_{i}\subseteq A_{j}}S_{i}}$ (notice that we may have: $S^{*}_{j}= \phi$). Let $S^{*}_{1},...,S^{*}_{\mid w^{*}\mid}$ denote the non empty sets $S^{*}_{j}$. Then $\chi^{*}=(S^{*}_{1},...,S^{*}_{\mid w^{*}\mid})$ is a smooth $(\frac{c_{0}}{\delta},\delta\lambda_{0},w^{*})$$-$structure of $T$ corresponding to $\mathcal{K}_{1}$ under $(\mathcal{K},\theta )$. Let $\epsilon_{1}$ be the $\epsilon$ from Lemma \ref{m} taken for $c= \frac{c_{0}}{\delta}$. Taking $\epsilon < \epsilon_{1}$ and since $\lambda_{0} < \frac{1}{\delta(2\delta)^{\delta +2}}$, then we can use Lemma \ref{m} and conclude taking $k=l$ that $\tilde{\mathcal{K}}_{1}$ is well-contained in $\chi^{*}$. Denote this well-contained copy of $\tilde{\mathcal{K}}_{1}$ by $G_{1}$. $\forall 1\leq i \leq p$, let $R_{i}^{*}=\displaystyle{\bigcap_{x\in V(G_{1})}}R_{i,x}$, where $R_{i} = \lbrace v \in \bigcup_{j=1}^{\mid w \mid}A_{j};$ $ \xi(v) = q_{i} \rbrace$. Let $ 1\leq i \leq p$. If $R_{i}= A_{j_{1}}$ for some $1\leq j_{1} \leq $ $\mid$$w$$\mid$ with $w(j_{1})=0$, then by Lemma \ref{g} $\mid$$R_{i}^{*}$$\mid$ $\geq (1-p\lambda_{0})\mid$$A_{j_{1}}$$\mid \geq (1-\delta\lambda_{0})\mid$$A_{j_{1}}$$\mid \geq \frac{\mid A_{j_{1}}\mid}{2} \geq \frac{\mid A_{j_{1}}\mid}{2\delta}$ since $\lambda_{0}\leq\frac{1}{2\delta}$. In this case we only rename the set $R_{i}^{*}$ by $R^{**}_{i}$. Let $ 1\leq i \leq p$. If $R_{i}\subseteq A_{j_{2}}$ for some $1\leq j_{2} \leq $ $\mid$$w$$\mid$ with $w(j_{2})=1$, then by Lemma \ref{g} $\mid$$R_{i}^{*}$$\mid$ $\geq (1-\delta^{2}\lambda_{0})\frac{\mid A_{j_{2}}\mid}{\delta}  \geq \frac{\mid A_{j_{2}}\mid}{2\delta}$ since $\lambda_{0}\leq\frac{1}{2\delta^{2}}$. In this case we select arbitrary $\lceil \frac{\mid A_{j_{2}}\mid}{2\delta}\rceil$ vertices from $R_{i}^{*}$ and we denote by $R_{i}^{**}$ the set of the selected $\lceil \frac{\mid A_{j_{2}}\mid}{2\delta}\rceil$ vertices. Now $\forall 1\leq j \leq $ $\mid$$w$$\mid$, let $M_{j}=  \displaystyle{\bigcup_{R^{**}_{i}\subseteq A_{j}}R^{**}_{i}}$ (notice that we may have: $M_{j}= \phi$). Let $M_{1},...,M_{\mid \overline{w}\mid}$ denote the non empty sets $M_{j}$. Also notice that for all $1\leq j\leq\mid$$w$$\mid$, $\mid$$M_{j}$$\mid$ $\geq \frac{\mid A_{s}\mid}{2\delta}$ for some $1\leq s\leq$ $\mid$$w$$\mid$. Then $\chi^{'} = (M_{1},...,M_{\mid \overline{w}\mid})$ form a smooth $(\frac{c_{0}}{2\delta},2\delta\lambda_{0},\overline{w})$$-$structure of $T$ corresponding to $\mathcal{K}_{2}$ under $(\mathcal{K},\theta )$ for an appropriate $\lbrace 0,1 \rbrace$$-$vector $\overline{w}$.  Let $\epsilon_{2}$ be the $\epsilon$ from Lemma \ref{o} taken for $c\geq \frac{c_{0}}{2\delta}$. Taking $\epsilon < \epsilon_{2}$ and since $\lambda_{0} < \frac{1}{(2\delta)^{\delta +3}}$, then we can use Lemma \ref{o} and conclude taking $k=m$ that $\tilde{\mathcal{K}}_{2}$ is well contained in $\chi^{'}$. Denote this copy well-contained of $\tilde{\mathcal{K}}_{2}$ by $G_{2}$. Now by merging $G_{1}$ and $G_{2}$ we get a well-contained copy of $\tilde{\mathcal{K}}$ in $\chi$. This completes the proof. $\blacksquare$ 
\end{proof}
\vspace{3mm}\\
From the previous lemma we get the following lemma:
\begin{lemma}
Let $\mathcal{N}$ be a regular super nebula under $\theta_{1}$ with $\mid$$\mathcal{N}$$\mid =\mu_{1}$ and let $\mathcal{G}$ be a regular $\Delta$galaxy under $\theta_{2}$ with $\mid$$\mathcal{G}$$\mid =\mu_{2}$. Let $\delta =\mu_{1}\mu_{2}$. Let $\mathcal{K}=\mathcal{N}\otimes\mathcal{G}$ under $\theta =(h_{1},...,h_{h})$ be a key tournament corresponding to $\mathcal{N}$ and $\mathcal{G}$ under $\theta_{1}$ and $\theta_{2}$ respectively. Let $0 < \lambda_{0} < \frac{1}{(2\delta)^{\delta+3}}$, $c_{0} > 0$ be constants, and let $w$ be a $\lbrace 0,1 \rbrace$$-$vector. Suppose that  $\chi = (S_{1},...,S_{\mid w \mid})$ is a smooth $(c_{0},\lambda_{0} ,w)$$-$structure of an $\epsilon$$-$critical tournament $T$ ($\epsilon$ is small enough) corresponding to $\mathcal{K}$ under the ordering $\theta$. Then 
\begin{itemize}
\item $T$ contains $\mathcal{N}$ or
\item $T$ contains $\mathcal{G}$.
\end{itemize}
\end{lemma}
\begin{proof}
Taking $\epsilon >0$ small enough, we conclude using the previous lemma that $\tilde{\mathcal{K}}$ is well-contained in $\chi$.  Denote by $G$ the well-contained copy of $\tilde{\mathcal{K}}$ in $\chi$. Let $\tilde{\theta}=(x_{1},...,x_{h})$ be the ordering of the vertices of $G$ according to their appearence in $\chi$ (i.e $\forall 1\leq i\leq h$, $\xi(x_{i})=i)$. Let $\mathcal{T}=T$$\mid$$V(G)$ and let $E=E(\mathcal{T})\backslash E(G)$. We have two cases: \\
Case $1$: For all $e\in E$, $e$ is a forward arc of $\mathcal{T}$ under $\tilde{\theta}$. Then $\mathcal{T}$ under $\tilde{\theta}$ is the key tournament $\mathcal{K}=\mathcal{N}\otimes\mathcal{G}$ under $\theta =(h_{1},...,h_{h})$ corresponding to $\mathcal{N}$ and $\mathcal{G}$ under $\theta_{1}$ and $\theta_{2}$ respectively (i.e $\tilde{\theta}$ is the super nebula ordering $\theta$ of $\mathcal{K}$). Then $\mathcal{T}$ under $\tilde{\theta}$ satisfies the following: 
 $V(G)=V(\mathcal{T})=\lbrace x_{1},...,x_{h}\rbrace$ is partitioned as follows: $V(\mathcal{T})= P\cup U_{1}\cup ... \cup U_{s}$, where $s=\mid$$B^{\mathcal{K},\theta}$$\mid$ $=\mid$$B^{\mathcal{T},\tilde{\theta}}$$\mid$  and $U_{i}$ corresponds to the bad triplet $t_{i}=(x_{i_{1}},x_{i_{2}},x_{i_{3}})\in B^{\mathcal{T},\tilde{\theta}}$ for $i=1,...,s$, such that:
 $\mathcal{T}$$\mid$$P$ is isomorphic to $\mathcal{N}$ and the restriction of $\tilde{\theta}$ to $P$ is the super nebula ordering $\theta_{1}$ of $\mathcal{N}$. 
  So $T$ contains $\mathcal{N}$ and we are done.\\
   The partition  of $V(G)$ that we get in case $1$ will be very useful and helpful in our analysis in case $2$ (i.e useful in identifying precisely the vertices of $G$, such that the subtournament of $T$ induced by these chosen vertices form a copy of $\mathcal{G}$). \\
 Case $2$: There exist $e_{r}=(x_{r_{j_{1}}},x_{r_{j_{2}}})\in E$ with $1\leq r\leq s$ and $j_{1},j_{2}\in\lbrace 1,2,3\rbrace$, such that $e_{r}$ is a backward arc of $\mathcal{T}$ under $\tilde{\theta}$. In this case $T$$\mid$$(U_{r}\cup\lbrace x_{r_{1}},x_{r_{2}},x_{r_{3}}\rbrace)$ is isomorphic to $\mathcal{G}$ and the restriction of $\tilde{\theta}$ to $U_{r}\cup\lbrace x_{r_{1}},x_{r_{2}},x_{r_{3}}\rbrace$ is the $\Delta$galaxy ordering $\theta_{2}$ of $\mathcal{G}$ (see the $5^{th}$ property of key tournaments in page $8$). So $T$ contains $\mathcal{G}$ and we are done. This completes the proof $\blacksquare$.
\end{proof}
\vspace{3mm}\\
Now we are ready to prove Theorem \ref{p}\vspace{3mm}\\
\noindent \sl {Proof of Theorem \ref{p}.} \upshape
Let $\mathcal{N}$ be a super nebula under $\theta_{1}$ with $\mid$$\mathcal{N}$$\mid =\mu_{1}$ and let $\mathcal{G}$ be a $\Delta$galaxy under $\theta_{2}$ with $\mid$$\mathcal{G}$$\mid =\mu_{2}$. Let $\delta =\mu_{1}\mu_{2}$ and let $0 < \lambda_{0} < \frac{1}{(2\delta)^{\delta+3}}$. We may assume that $\mathcal{N}$ is a regular super nebula and $\mathcal{G}$ is a regular $\Delta$galaxy since every super nebula is a subtournament of a regular super nebula and every $\Delta$galaxy is a subtournament of a regular $\Delta$galaxy.  Let $\mathcal{K}=\mathcal{N}\otimes\mathcal{G}$ under $\theta $ be a key tournament corresponding to $\mathcal{N}$ and $\mathcal{G}$ under $\theta_{1}$ and $\theta_{2}$ respectively.
 Let $\epsilon > 0$ be small enough. Assume that $\lbrace \mathcal{N},\mathcal{G}\rbrace$ does not satisfy $EHC$, then there exists an $\lbrace \mathcal{N},\mathcal{G}\rbrace$$-$free $\epsilon$$-$critical tournament $T$. By Lemma \ref{e}, $\mid$$T$$\mid$ is large enough. By Theorem \ref{i}, $T$ contains a smooth $(c_{0},\lambda_{0},w)$$-$structure $(S_{1},...,S_{\mid w \mid})$ corresponding to $\mathcal{K}$ under $\theta$ for some $c_{0} >0$ and appropriate $\lbrace 0,1 \rbrace$$-$vector $w$. Then by the previous lemma, $T$ contains $\mathcal{N}$ or $T$ contains $\mathcal{G}$, a contradiction. $\blacksquare$   
\section{Central triangular galaxies and $K_{6}$} 
\subsection{Definitions} \label{d}
Let $ \beta = (v_{1},...,v_{f}) $ be an ordering of the vertex set $V(T)$ of an $ f- $vertex tournament $T$. A $K_{6}=\lbrace v_{i_{1}},...,v_{i_{6}}\rbrace$ of $T$ under $\beta$ (where $i_{1}<...<i_{6}$) is the subtournament of $T$ induced by $\lbrace v_{i_{1}},...,v_{i_{6}}\rbrace$ such that $T$$\mid$$\lbrace v_{i_{1}},...,v_{i_{6}}\rbrace$ is the tournament $K_{6}$, $ (v_{i_{1}},...,v_{i_{6}})$ is the canonical ordering of $K_{6}$, $i_{2}=i_{1}+1$, $i_{4}=i_{3}+1$, and $i_{6}=i_{5}+1$. We call $v_{i_{1}}$ and $v_{i_{6}}$ the \textit{centers of} $K_{6}$.\\
Let $K_{6}^{i}=\lbrace v_{i_{1}},...,v_{i_{6}}\rbrace$ be a $K_{6}$ of $T$ under $\beta$. Define $operation_{K_{6},\theta}$ by  deleting the vertices $v_{i_{2}},v_{i_{3}},v_{i_{5}} $ and reversing the orientation of the arc $(v_{i_{4}},v_{i_{6}}) $. This $K_{6},\beta$ in $operation_{K_{6},\beta}$ is because this operation is applied for the tournaments $K_{6}$ of $T$ under $\beta$. \\
A tournament $T$ is a $GK_{6}$ if there exist an ordering  $ \beta = (v_{1},...,v_{f}) $ of its vertices such that $V(T)$ is the disjoint union of $V(K^{1}_{6}),...,V(K^{l}_{6}),X$, and such that $K^{1}_{6},...,K^{l}_{6}$ are the $K_{6}$ tournaments of $T$ under $\beta$, $T$$\mid$$X$ is a regular galaxy under the restriction of $\beta$ to $X$, and no center of a $K_{6}$ of $T$ under $\beta$ is between leaves of a star of $T$ under $\beta$. In this case we also say that $T$ is a $GK_{6}$ \textit{under} $\beta$. Obviously notice that every $GK_{6}$ under $\beta$ is a super nebula under $\beta$.\vspace{3mm}

\begin{figure}[h]
	\centering
	\includegraphics[width=0.7\linewidth]{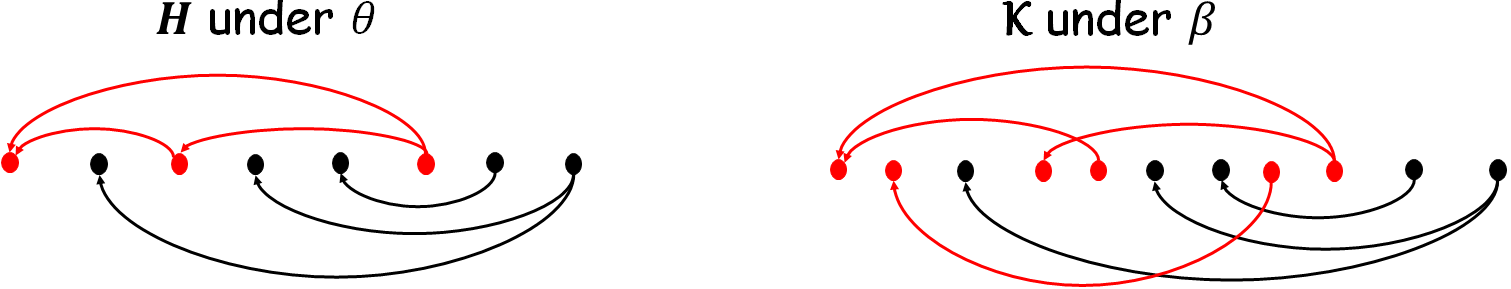}
	\caption{Central triangular galaxy $H$ under $\theta$ and its corresponding key tournament $\mathcal{K}$ under $\beta$. All the non drawn arcs are forward.}
	\label{fig:keyk6tour}
\end{figure}

Let $H$ be a regular central triangular galaxy with $\mid$$ H $$\mid$ $= h$ and let $(u_{1},...,u_{h})$ be a central triangular galaxy ordering of $H$. Denote this ordering by $\theta$. Let $\Delta_{1},...,\Delta_{l}$ be the triangles of $H$ under $\theta$ and let $Q_{1},...,Q_{m}$ be the frontier stars of $H$ under $\theta $.  A $GK_{6}$ tournament $\mathcal{K}$ under $\beta$ is a \textit{key tournament corresponding to $H$ under $\theta $} if the tournament obtained from $\mathcal{K}$ under $\beta$ after performing $operation_{K_{6},\beta}$ to every $K_{6}$ of $\mathcal{K}$ under $\beta$ is the tournament $H$ and the obtained ordering is the central triangular galaxy ordering $\theta$ of $H$ (see Figure \ref{fig:keyk6tour}).

Let $H$ be a regular central triangular galaxy under $\theta$ with $\mid$$ H $$\mid$ $= h$. Let $\Delta_{1},...,\Delta_{l}$ be the triangles of $H$ under $\theta$ and let $Q_{1},...,Q_{m}$ be the frontier stars of $H$ under $\theta $. Assume without loss of generality that $m=l$. Let $\mathcal{K}$ under $\beta = (v_{1},...,v_{f})$ be the key tournament corresponding to $H$ under $\theta $. Let $K^{1}_{6},...,K^{l}_{6}$ be the $K_{6}$ tournaments of $\mathcal{K}$ under $\beta$. Let $X=V(\mathcal{K})\backslash \bigcup_{i=1}^{l}V(K_{6}^{i})$ and let $Q_{1},...,Q_{l}$ be the stars of $\mathcal{K}$$\mid$$X$ under $\beta_{X}$ ($\beta_{X}$ is the restriction of $\beta$ to $X$). Notice that $f= h+3l$. For $k \in \lbrace 0,...,l \rbrace$ define $\mathcal{K}^{k} = \mathcal{K}$$\mid$$\bigcup_{j=1}^{k}(V(K^{j}_{6})\cup V(Q_{j}))$ where $\mathcal{K}^{l} = \mathcal{K}$, and $\mathcal{K}^{0}$ is the empty tournament. For $k \in \lbrace 1,...,l \rbrace$ let $\beta_{k}=(v_{k_{1}},...,v_{k_{q_{k}}})$ be the restriction of $\beta$ to $V(\mathcal{K}^{k})$.  Recall that $s^{\mathcal{K},\beta}$ is a $\lbrace 0,1 \rbrace$$-$vector such that $s^{\mathcal{K},\beta}(i) = 0$ if and only if $v_{i}\in \mathcal{C}$, where $\mathcal{C}$ is the set of all chosen centers of the stars of $\mathcal{K}$ under $\beta$ and the centers of the super middle $2$-nebulas of $\mathcal{K}$ under $\beta$. Let $s^{\mathcal{K},\beta}_{\mathcal{K}^{k}}$ be the restriction of $s^{\mathcal{K},\beta}$ to the $0's$ and $1's$ corresponding to $V(\mathcal{K}^{k})$ (notice that $s^{\mathcal{K},\beta}_{\mathcal{K}^{k}}= s^{\mathcal{K}^{k},\beta_{k}}$) and let $^{c}s^{\mathcal{K},\beta}_{\mathcal{K}^{k}}$ be the vector obtained from $s^{\mathcal{K},\beta}_{\mathcal{K}^{k}}$ by replacing every subsequence of consecutive $1's$ corresponding to the same entry of $s^{\mathcal{K},\beta}_{c}$ by single $1$ (see page $7$ for the definition of $s^{\mathcal{K},\beta}_{c}$). We say that a smooth $(c,\lambda ,w)$$-$structure of a tournament $T$ \textit{corresponds}  \textit{to $\mathcal{K}^{k}$ under $(\mathcal{K},\beta )$} if $w = ^{c}s^{\mathcal{K},\beta}_{\mathcal{K}^{k}}$. Notice that $s^{\mathcal{K},\beta}_{\mathcal{K}^{l}}=s^{\mathcal{K},\beta}$ and $^{c}s^{\mathcal{K},\beta}_{\mathcal{K}^{l}}=s^{\mathcal{K},\beta}_{c}$.\\
Let $\nu =^{c}s^{\mathcal{K},\beta}_{\mathcal{K}^{k}}$. Let $\delta^{\nu}:$ $\lbrace j: \nu_{j} = 1 \rbrace \rightarrow \mathbb{N}$ be a function that assigns to every nonzero entry of $\nu$ the number of consecutive $1'$s of $s^{\mathcal{K},\beta}_{\mathcal{K}^{k}}$ replaced by that entry of $\nu$.\\
Fix $k \in \lbrace 0,...,l \rbrace$. Let $\tilde{\mathcal{K}}^{k}= \tilde{\mathcal{K}}$$\mid$$V(\mathcal{K}^{k})$, where $\tilde{\mathcal{K}}$ is the mutant super nubula obtained from $\mathcal{K}$ under $\theta$ (see page $7$). Let $(S_{1},...,S_{\mid w \mid})$ be a smooth $(c,\lambda ,w)$$-$structure of a tournament $T$ corresponding to $\mathcal{K}^{k}$ under $(\mathcal{K},\beta )$. Let $i_{r}$ be such that $w(i_{r}) = 1$. Assume that $S_{i_{r}} = \lbrace s^{1}_{i_{r}},...,s_{i_{r}}^{\mid S_{i_{r}} \mid} \rbrace$ and $(s^{1}_{i_{r}},...,s_{i_{r}}^{\mid S_{i_{r}} \mid})$ is a transitive ordering. Write $m(i_{r}) = \lfloor\frac{\mid S_{i_{r}} \mid}{\delta^{w}(i_{r})}\rfloor$. Denote $S^{j}_{i_{r}} = \lbrace s^{(j-1)m(i_{r})+1}_{i_{r}},...,s_{i_{r}}^{jm(i_{r})} \rbrace$ for $j \in \lbrace 1,...,\delta^{w}(i_{r}) \rbrace$. For every $v \in S^{j}_{i_{r}}$ denote $\xi(v) = (\mid$$\lbrace k < i_{r}: w(k) = 0 \rbrace$$\mid$ $+$ $\displaystyle{\sum_{k < i_{r}: w(k) = 1}\delta^{w}(k) })$ $+$ $j$. For every $v \in S_{i_{r}}$ such that $w(i_{r}) = 0$ denote $\xi(v) = (\mid$$\lbrace k < i_{r}: w(k) = 0 \rbrace$$\mid$ $+$ $\displaystyle{\sum_{k < i_{r}: w(k) = 1}\delta^{w}(k) })$ $+$ $1$. We say that $\tilde{\mathcal{K}}^{k}$ is \textit{well-contained in} $(S_{1},...,S_{\mid w \mid})$ that corresponds to $\mathcal{K}^{k}$ under $(\mathcal{K},\beta )$ if there is an injective homomorphism $f$ of $\tilde{\mathcal{K}}^{k}$ into $T$$\mid$$\bigcup_{i = 1}^{\mid w \mid}S_{i}$ such that $\xi(f(v_{k_{j}})) = j$ for every $j \in \lbrace 1,...,q_{k} \rbrace$.

\subsection{Proof of Theorem \ref{t}}
We start by the following technical lemma:
\begin{lemma}
Let $H$ be a regular central triangular galaxy under $\theta$ with $\mid$$H$$\mid$ $= h$. Let $\Delta_{1},...,\Delta_{l}$ be the triangles of $H$ under $\theta$ and let $Q_{1},...,Q_{l}$ be the stars of $H$ under $\theta $. Let $\mathcal{K}$ under $\beta$ be the key tournament corresponding to $H$ under $\theta $ ($\mid$$\mathcal{K}$$\mid$ $= h+3l$). Let $K^{1}_{6},...,K^{l}_{6}$ be the $K_{6}$ tournaments of $\mathcal{K}$ under $\beta$ and let $Q_{1},...,Q_{l}$ be the stars of $\mathcal{K}$$\mid$$X$ under $\beta_{X}$, where $X=V(\mathcal{K})\backslash \bigcup_{i=1}^{l}V(K_{6}^{i})$ and $\beta_{X}$ is the restriction of $\beta$ to $X$. Let $0 < \lambda < \frac{1}{(4h)^{h+4}}$, $c > 0$ be constants, and $w$ be a $\lbrace 0,1 \rbrace$$-$vector. Fix $k \in \lbrace 0,...,l \rbrace$ and let $\widehat{\lambda} = (4h)^{l-k}\lambda$ and $\widehat{c} = \frac{c}{(4h)^{l-k}}$. There exist $ \epsilon_{k} > 0$ such that $\forall 0 < \epsilon < \epsilon_{k}$, for every $\epsilon$$-$critical tournament $T$ with $\mid$$T$$\mid$ $= n$ containing $\chi = (S_{1},...,S_{\mid w \mid})$ as a smooth $(\widehat{c},\widehat{\lambda},w)$$-$structure corresponding to $\mathcal{K}^{k}$ under $(\mathcal{K},\beta )$, we have $\tilde{\mathcal{K}}^{k}$ is well-contained in $\chi$.    
\end{lemma}
\begin{proof}
The proof is by induction on $k$. For $k=0$ the statement is obvious since $\tilde{\mathcal{K}}^{0}$ is the empty digraph. Suppose that $\chi = (S_{1},...,S_{\mid w \mid})$ is a smooth $(\widehat{c},\widehat{\lambda},w)$$-$structure in $T$ corresponding to $\mathcal{K}^{k}$ under $(\mathcal{K},\beta )$ with $\beta= (h_{1},...,h_{\mid\mathcal{K}\mid})$ and $\mid$$\mathcal{K}$$\mid$ $= h+3l$. Let $\beta_{k}=(h_{k_{1}},...,h_{k_{p}})$ be the restriction of $\beta$ to $V(\mathcal{K}^{k})$. Let $K_{6}^{k} = \lbrace h_{k_{s_{1}}},...,h_{k_{s_{6}}} \rbrace$. Assume without loss of generality that the star $\lbrace h_{k_{s_{2}}},h_{k_{s_{5}}} \rbrace$ is considered as a left star of $\mathcal{K}$ under $\beta$.  Let $h_{k_{p_{0}}}$ be the center of $Q_{k}$ and $h_{k_{p_{1}}},...,h_{k_{p_{q}}}$ be its leafs for some integer $q>0$. Let $D_{i} = \lbrace v \in \bigcup_{j=1}^{\mid w \mid}S_{j};$ $ \xi(v) = s_{i} \rbrace$ for $i=1,...,6$. Then $\exists x_{1},x_{2},x_{3},z_{1},z_{2} \in \lbrace 1,...,\mid$$w$$\mid \rbrace$ with $x_{1}<x_{2}<z_{1}\leq z_{2}<x_{3}$, $w(x_{1})=w(x_{2})=w(x_{3})=0$, and $w(z_{1})=w(z_{2})=1$, such that $D_{1} = S_{x_{1}}$, $D_{2} = S_{x_{2}}$, $D_{i} \subseteq S_{z_{1}}$ for $i=3,4$,  $D_{5} \subseteq S_{z_{2}}$, and $D_{6} = S_{x_{3}}$. $\forall 0 \leq i \leq q$, let $R_{i} = \lbrace v \in \bigcup_{j=1}^{\mid w \mid}S_{j};$ $ \xi(v) = p_{i} \rbrace$. Then $\exists x_{4} \in \lbrace 1,...,\mid$$w$$\mid \rbrace \backslash \lbrace x_{1},x_{2},x_{3},z_{1},z_{2} \rbrace$ with $w(x_{4})=0$, and $\exists z_{3} \in \lbrace 1,...,\mid$$w$$\mid \rbrace \backslash \lbrace x_{1},x_{2},x_{3},x_{4} \rbrace$ with $w(z_{3})=1$, such that $R_{0} = S_{x_{4}}$ and $\forall 1 \leq i \leq q$, $R_{i} \subseteq S_{z_{3}}$. Since we can assume that $\epsilon < min \lbrace log_{\frac{\widehat{c}}{4}}(1-\frac{\widehat{c}}{2h}),log_{\frac{\widehat{c}}{4}}(\frac{1}{2})\rbrace$, then by Lemma \ref{s} there exists vertices $d_{1},d_{3},d_{4},d_{6}$ such that $d_{i} \in D_{i}$ for $i=1,3,4,6$ and $d_{1}\leftarrow d_{4}$ and $\lbrace d_{1},d_{3}\rbrace \leftarrow d_{6}$. Also notice that $d_{3}\rightarrow d_{4}$. One of the following holds:\\
$\bullet$ $z_{1}< z_{2}$ and $z_{3}\notin \lbrace z_{1},z_{2}\rbrace$, or\\
$\bullet$ $z_{1}< z_{2}$ and $z_{3}=z_{1}$, or\\
$\bullet$ $z_{1}< z_{2}$ and $z_{3}=z_{2}$, or\\
$\bullet$ $z_{1}= z_{2}$ and $z_{3}\neq z_{1}$, or\\
$\bullet$ $z_{1}= z_{2}=z_{3}$.\\
Assume that $z_{1}< z_{2}$ and $z_{3}\notin \lbrace z_{1},z_{2}\rbrace$. Else, the argument is similar and we omit it.\\ Let $D_{2}^{*} = \lbrace d_{2}\in D_{2}; d_{1}\rightarrow d_{2} \rightarrow \lbrace d_{3},d_{4},d_{6}\rbrace \rbrace$ and $D_{5}^{*} = \lbrace d_{5}\in D_{5}; \lbrace d_{1},d_{3},d_{4}\rbrace\rightarrow d_{5} \rightarrow d_{6} \rbrace$. Then by Lemma \ref{g}, $\mid$$D_{2}^{*}$$\mid$ $\geq (1-4\widehat{\lambda})\widehat{c}n\geq \frac{\widehat{c}}{2}n$ since $\widehat{\lambda} \leq \frac{1}{8}$, and $\mid$$D_{5}^{*}$$\mid$ $\geq \frac{1-8h\widehat{\lambda}}{2h}\widehat{c}tr(T)\geq \frac{\widehat{c}}{4h}tr(T)$ since $\lambda \leq \frac{1}{16h}$. 
 Since we can assume that $\epsilon < log_{\frac{\widehat{c}}{2}}(1-\frac{\widehat{c}}{4h})$, then Lemma \ref{f} implies that there exist vertices $d_{2}\in D_{2}^{*}$ and $d_{5}\in D_{5}^{*}$ such that $d_{2} \leftarrow d_{5}$. Then $T$$\mid$$\lbrace d_{1},...,d_{6}\rbrace$ contains a copy of $\tilde{\mathcal{K}}^{k}$$\mid$$V(K_{6}^{k})$. Denote this copy by $W$.   
$\forall 0 \leq i \leq q$, let $R_{i}^{*} = \bigcap_{x\in V(W)}R_{i,x}$.
  Then by Lemma \ref{g}, $\mid$$R_{0}^{*}$$\mid$ $\geq (1-6\widehat{\lambda})\mid$$R_{0}$$\mid$ $\geq \frac{\mid R_{0}\mid}{2}$ $\geq\frac{\widehat{c}}{2}n$ since $\widehat{\lambda} \leq \frac{1}{12}$, and $\forall 1 \leq i \leq q$, $\mid$$ R_{i}^{*}$$\mid$ $ \geq \frac{1-12h\widehat{\lambda}}{2h}\mid$$ S_{z_{3}}$$ \mid$ $\geq \frac{\widehat{c}}{4h}tr(T)$ since $\widehat{\lambda} \leq \frac{1}{24h}$. Since we can assume that $\epsilon < log_{\frac{\widehat{c}}{4h}}(1-\frac{\widehat{c}}{4h})$, then by Lemma \ref{r} there exists vertices $r_{0},r_{1},...,r_{q}$ such that $r_{i} \in R_{i}^{*}$ for $i=0,1,...,q$ and \\
$\ast$ $r_{1},...,r_{q}$ are all adjacent from $r_{0}$ if $x_{4}>z_{3}$.\\
$\ast$ $r_{1},...,r_{q}$ are all adjacent to $r_{0}$ if $x_{4}<z_{3}$.\\
So $T$$\mid$$\lbrace d_{1},...,d_{6},r_{0},r_{1},...,r_{q} \rbrace$ contains a copy of $\tilde{\mathcal{K}}^{k}$$\mid$$(V(K_{6}^{k})\cup V(Q_{k}))$. Denote this copy by $Y$.\\
$\forall i \in \lbrace 1,...,\mid$$w$$\mid \rbrace \backslash \lbrace x_{1},...,x_{4},z_{1},z_{2},z_{3} \rbrace$, let $S_{i}^{*} = \displaystyle{\bigcap_{x\in V(Y)}S_{i,x}}$.  Then by Lemma \ref{g}, $\mid$$S_{i}^{*}$$\mid$ $\geq (1-\mid$$Y$$\mid\widehat{\lambda})\mid$$S_{i}$$\mid$ $\geq (1-2h\widehat{\lambda})\mid$$S_{i}$$\mid$ $\geq \frac{1}{4h}\mid$$S_{i}$$\mid$ since $\widehat{\lambda} \leq \frac{4h-1}{8h^{2}}$.
Write $\mathcal{H} = \lbrace 1,...,p \rbrace \backslash \lbrace p_{0},...,p_{q},s_{1},...,s_{6} \rbrace$. Let $Z_{i}=\lbrace v\in V(Y): v\in S_{z_{i}}\rbrace$ for $i=1,2,3$. $\forall 1\leq i\leq 3$, if $\lbrace v\in S_{z_{i}}: \xi(v) \in \mathcal{H} \rbrace \neq \phi$, then define $J_{z_{i}} = \lbrace \eta \in \mathcal{H}: \exists v \in S_{z_{i}}$ and $\xi(v)= \eta \rbrace$. Now $\forall \eta \in J_{z_{i}}$, let $S_{z_{i}}^{*\eta}= \lbrace v \in S_{z_{i}}: \xi(v)=\eta$ and $v \in \displaystyle{\bigcap_{x\in V(Y)\backslash Z_{i}}S_{z_{i},x}} \rbrace$. Then by Lemma \ref{g}, $\forall \eta \in J_{z_{i}}$, we have $\mid$$S_{z_{i}}^{*\eta}$$\mid$ $\geq \frac{1-2h^{2}\widehat{\lambda}}{2h}\mid $$S_{z_{i}}$$\mid $ $\geq \frac{\mid S_{z_{i}}\mid}{4h}$ since $\widehat{\lambda} \leq \frac{1}{4h^{2}}$. Now $\forall \eta \in J_{z_{i}}$, select arbitrary $\lceil \frac{\mid S_{z_{i}}\mid}{4h}\rceil$ vertices of $S_{z_{i}}^{*\eta}$ and denote the union of these $\mid$$J_{z_{i}}$$\mid$ sets by $S_{z_{i}}^{*}$. 
So we have defined some number of sets. Denote by $t$ the number of these defined sets and by $S_{1}^{*},...,S^{*}_{t}$ these  sets. We have $\mid$$S_{i}^{*}$$\mid$ $\geq \frac{S_{i}}{4h}$ for every defined set. Now Lemma \ref{b} implies that $\chi^{*}=(S_{1}^{*},...,S^{*}_{t})$ form a smooth $(\frac{\widehat{c}}{4h},4h\widehat{\lambda},w^{*})$$-$structure of $T$ corresponding to $\mathcal{K}^{k-1}$ under $(\mathcal{K},\beta )$, where $\frac{\widehat{c}}{4h}= \frac{c}{(4h)^{l-(k-1)}}, 4h\widehat{\lambda}=(4h)^{l-(k-1)}\lambda$, and $w^{*}$ is an appropriate $\lbrace 0,1 \rbrace$$-$vector.
Now take $\epsilon_{k} < min \lbrace \epsilon_{k-1}, log_{\frac{\widehat{c}}{4}}(1-\frac{\widehat{c}}{2h}),log_{\frac{\widehat{c}}{4}}(\frac{1}{2}),log_{\frac{\widehat{c}}{4h}}(1-\frac{\widehat{c}}{4h}) \rbrace$. So by induction hypothesis $\tilde{\mathcal{K}}^{k-1}$ is well-contained in $\chi^{*}$. Now by merging the well-contained copy of $\tilde{\mathcal{K}}^{k-1}$ and $Y$ we get a well-contained copy of $\tilde{\mathcal{K}}^{k}$. $\blacksquare$\vspace{3mm}\\ 
    \end{proof}
    From the above lemma we get the following lemma: 
    \begin{lemma}
Let $H$ be a regular central triangular galaxy under $\theta$ with $\mid$$H$$\mid$ $= h$. Let $\mathcal{K}$ under $\beta$ be the key tournament corresponding to $H$ under $\theta$. Let $0 < \lambda < \frac{1}{(4h)^{h+4}}$, $c > 0$ be constants, and let $w$ be a $\lbrace 0,1 \rbrace$$-$vector. Suppose that  $\chi = (S_{1},...,S_{\mid w \mid})$ is a smooth $(c,\lambda ,w)$$-$structure of an $\epsilon$$-$critical tournament $T$ ($\epsilon$ is small enough) corresponding to $\mathcal{K}$ under the ordering $\beta$. Then  $T$ contains $H$ or $T$ contains $K_{6}$.
\end{lemma}
\begin{proof}
Let $K^{1}_{6},...,K^{l}_{6}$ be the $K_{6}$ tournaments of $\mathcal{K}$ under $\beta$ and let $Q_{1},...,Q_{l}$ be the stars of $\mathcal{K}$$\mid$$X$ under $\beta_{X}$. Taking $\epsilon >0$ small enough and $k=l$, we conclude using the previous lemma that $\tilde{\mathcal{K}}$ is well-contained in $\chi$.  Denote by $G$ the well-contained copy of $\tilde{\mathcal{K}}$ in $\chi$. $\forall 1\leq i\leq l$, let $\mathcal{D}_{i}=\lbrace d_{1}^{i},...,d_{6}^{i}\rbrace$ be the copy of $\tilde{\mathcal{K}}$$\mid$$V(K_{6}^{i})$ in $\chi$ and let $\mathcal{Q}_{i}$ be the copy of $\tilde{\mathcal{K}}$$\mid$$V(Q_{i})$ in $\chi$. Let $\tilde{\theta}$ be the ordering of the vertices of $G$ according to their appearence in $\chi$. Notice that $\forall 1\leq i\leq l$, we don't know the orientation of the edges $d_{4}^{i}d_{6}^{i}$ and $d_{1}^{i}d_{3}^{i}$.\\
Assume first that $\forall 1\leq i\leq l$, at least one of the following holds:\\
$\bullet$ $d_{1}^{i}\leftarrow d_{3}^{i}$\\
$\bullet$ $d_{4}^{i}\leftarrow d_{6}^{i}$\\   
Let $j_{i}\in\lbrace 3,4\rbrace$, such that $d_{1}^{i}\leftarrow d_{j_{i}}^{i}\leftarrow d_{6}^{i}$ for $i=1,...,l$. But then  the restriction of $\tilde{\theta}$ to $\bigcup_{i=1}^{l}(V(\mathcal{Q}_{i})\cup \lbrace d_{1}^{i}, d_{j_{i}}^{i}, d_{6}^{i}\rbrace)$ is the central triangular galaxy ordering $\theta$ of $H$. So $T$ contains $H$, and we are done. \\ Otherwise there exist $i\in \lbrace 1,...,l\rbrace$ such that  $d_{1}^{i}\rightarrow d_{3}^{i}$ and $d_{4}^{i}\rightarrow d_{6}^{i}$. But then $(d_{1}^{i},d_{2}^{i},d_{3}^{i},d_{4}^{i},d_{5}^{i}, d_{6}^{i})$ is the canonical ordering of $K_{6}$. So $T$ containes $K_{6}$. This completes the proof. $\blacksquare$\vspace{3mm}\\
\end{proof}
Now we are ready to prove Theorem \ref{t}:\vspace{3mm}\\
\noindent \sl {Proof of Theorem \ref{t}.} \upshape
Let $H$ be a regular central triangular galaxy under $\theta$ with $\mid$$H$$\mid$ $= h$. We may assume that $H$ is a regular central triangular galaxy since every central triangular galaxy is a subtournament of a regular central triangular galaxy. Let $\mathcal{K}$ under $\beta$ be the key tournament corresponding to $H$ under $\theta$. Let $\epsilon > 0$ be small enough and let $0 < \lambda < \frac{1}{(4h)^{h+4}}$ be constants.   
  Assume that $\lbrace H,K_{6}\rbrace$ does not satisfy $EHC$, then there exists an $\lbrace H,K_{6}\rbrace$$-$free $\epsilon$$-$critical tournament $T$. By Lemma \ref{e}, $\mid$$T$$\mid$ is large enough. By Theorem \ref{i}, $T$ contains a smooth $(c,\lambda ,w)$$-$structure $(S_{1},...,S_{\mid w \mid})$ corresponding to $\mathcal{K}$ under under $(\mathcal{K},\beta )$, for some $c >0$ and appropriate $\lbrace 0,1 \rbrace$$-$vector $w$. Then by the previous lemma, $T$ contains $H$ or $T$ contains $K_{6}$, a contradiction. $\blacksquare$ 
 
\section{Extension of the results} \label{k}
\subsection{Extension of Theorem \ref{p}}
Let $S = \lbrace u_{1},u_{2},...,u_{p}\rbrace$ be a middle star and let $u_{r}$ be the center of $S$ with $2\leq r\leq p-1$ (note that $(u_{1},u_{2},...,u_{p})$ is its star ordering). If $r=2$ then $S$ is called \textit{$1$-left middle star} and if $r= p-1$ then $S$ is called \textit{$1$-right middle star}.
\begin{theorem} \label{q}
Let $\mathcal{N}$ be a super nebula under $\theta$ and besides for every star $Q_{i}$ of $\mathcal{N}$ under $\theta$, $Q_{i}$ is a $1$-right middle star or a right star, and all the super $2$-nebulas of $\mathcal{N}$ under $\theta$ are left super $2$-nebulas. Let $\mathcal{G}$ be a left triangular galaxy under $\alpha$ and besides $\mathcal{G}$ has only one triangle under $\alpha$. Then $\lbrace \mathcal{N},\mathcal{G}\rbrace$ satisfy EHC. 
\end{theorem}
\begin{theorem} \label{aa}
Let $\mathcal{N}$ be a super nebula under $\theta$ and besides for every star $Q_{i}$ of $\mathcal{N}$ under $\theta$, $Q_{i}$ is a $1$-left middle star star or a left star, and all the super $2$-nebulas of $\mathcal{N}$ under $\theta$ are right super $2$-nebulas. Let $\mathcal{G}$ be a right triangular galaxy under $\alpha$ and besides $\mathcal{G}$ has only one triangle under $\alpha$. Then $\lbrace \mathcal{N},\mathcal{G}\rbrace$ satisfy EHC. 
\end{theorem}
\begin{theorem}\label{bb} 
Let $\mathcal{N}$ be a super nebula under $\theta$ and besides all the stars of $\mathcal{N}$ under $\theta$ are frontier stars, and for every super $2$-nebula $\Sigma_{i}$ of $\mathcal{N}$ under $\theta$, $\Sigma_{i}$ is a middle super $2$-nebula. Let $\mathcal{G}$ be a central triangular galaxy under $\alpha$ and besides $\mathcal{G}$ has only one triangle under $\alpha$. Then $\lbrace \mathcal{N},\mathcal{G}\rbrace$ satisfy EHC. 
\end{theorem}
We say that a tournament $H$ is a \textit{super $\Delta$galaxy under} $\theta$ if it is a super triangular galaxy under $\theta$ and besides $H$ has only one triangle under $\theta$. We say that a tournament $H$ is an \textit{$LR$-$\Delta$galaxy under} $\theta$ if it is a super $\Delta$galaxy under $\theta$ and besides the vertices of the triangle $\Delta$ of $H$ under $\theta$ that are allowed to be in the ordering $\theta$ between leaves of a star of $H$ under $\theta$ are only the exteriors of $\Delta$. We say that a tournament $H$ is a \textit{$CR$-$\Delta$galaxy} (resp. \textit{$CL$-$\Delta$galaxy}) \textit{under} $\theta$ if it is a super $\Delta$galaxy under $\theta$ and besides the vertices of the triangle $\Delta$ of $H$ under $\theta$ that are allowed to be in the ordering $\theta$ between leaves of a star of $H$ under $\theta$ are the right (resp. left) exterior and the center of $\Delta$, such that: if $Q_{i}$ and $Q_{j}$ are frontier stars of $H$ under $\theta$, such that the center of $\Delta$ is in the ordering between the leaves of $Q_{i}$ and the right (resp. left) exterior of $\Delta$ is in the ordering between the leaves of $Q_{j}$, then no leaf of $Q_{i}$ is between the leaves of $Q_{j}$ under $\theta$ and no leaf of $Q_{j}$ is between the leaves of $Q_{i}$ under $\theta$. And if the center of $\Delta$ is between the leaves of $Q_{i}$ for some star $Q_{i}$ of $H$ under $\theta$, then the right (resp. left) exterior of $\Delta$ is not between the leaves of $Q_{i}$ under $\theta$.\\
 We say that a tournament $H$ is a \textit{super left nebula} (resp. \textit{super right nebula}) \textit{under} $\theta$ if it is a nebula under $\theta$ and besides all the stars of $H$ under $\theta$ are left stars (resp. right stars).
\begin{theorem}\label{u}
If $H_{1}$ and $H_{2}$ are: a central nebula and a  $LR$-$\Delta$galaxy, or: a super left nebula and a  $CR$-$\Delta$galaxy, or: a super right nebula and a  $CL$-$\Delta$galaxy, then $\lbrace H_{1},H_{2}\rbrace$ satisfies the Erd\"{o}s-Hajnal Conjecture.  
\end{theorem}
We omit the proof of Theorems \ref{q}, \ref{aa}, \ref{bb}, \ref{u} because they have completely the same proof of Theorem \ref{p}.
\subsection{Generalization of Theorem \ref{t}}
 Let $H$ be a tournament such that there exists an ordering $\theta$ of its vertices such that $V(H)$ is the disjoint union of $V(\Sigma),V(Q_{1}),...,V(Q_{m})$ where $Q_{1},...,Q_{m}$ are the frontier stars of $H$ under $\theta$, $\Sigma$ is the  super $2$-nebula of $H$ under $\theta$, no center of a star is between leaves of $\Sigma$ under $\theta$, no center of $\Sigma$ is between leaves of a star of $H$ under $\theta$, and no center of a star appears in the ordering $\theta$ between leaves of another star. In this case $H$ is  called a \textit{$\Sigma$-galaxy under} $\theta$ and $\theta$ is called a \textit{$\Sigma$-galaxy ordering of $H$}. If $\Sigma$ is a super middle $2$-nebula (resp. super left $2$-nebula) (resp. super right $2$-nebula) then $H$ is called a \textit{middle $\Sigma$-galaxy} (resp. \textit{left $\Sigma$-galaxy}) (resp. \textit{right $\Sigma$-galaxy}).  
 Obviously one can notice that $K_{6}$ is a middle $\Sigma$-galaxy and its canonical ordering is its $\Sigma$-galaxy ordering. Also notice that every $\Sigma$-galaxy is a super nebula. The following theorem is a generalization of Theorem \ref{t}:
\begin{theorem} \label{l}
If $H_{1}$ and $H_{2}$ are: a middle $\Sigma$-galaxy and a  central triangular galaxy, or: a left $\Sigma$-galaxy and a  left triangular galaxy, or: a right $\Sigma$-galaxy and a  right triangular galaxy, then $\lbrace H_{1},H_{2}\rbrace$ satisfies the Erd\"{o}s-Hajnal Conjecture.  
\end{theorem}  
We omit the proof of Theorem \ref{l}, since it is completely analogous to the proof of Theorem \ref{t}. The proof uses the notion of key tournaments $\mathcal{K}$ under $\theta$ corresponding to $H_{1}$ under its $\Sigma$-galaxy ordering and  $H_{2}$ under its super triangular galaxy ordering. The problem we face is that when looking for a well-contained copy of $H_{1}$ (resp. $H_{2}$) in an appropriate smooth $(c,\lambda ,w)$-structure, there are a group of arcs that we know nothing about their orientation. This is the place where we need to use key tournaments constructed depending on both $H_{1}$ and $H_{2}$ (we construct it following the same principle in Section \ref{d}).  We first find this mutant key tournament $\tilde{\mathcal{K}}$ as a well-contained copy in a smooth $(c,\lambda ,w)$-structure corresponding to $\mathcal{K}$ under $\theta$. Then we extract $H_{1}$ or $H_{2}$ depending on the orientation of the arcs where the problem is faced.


\begin{thebibliography}{99}
\bibitem{fdo}N. Alon, J. Pach, J. Solymosi, Ramsey$ - $type theorems with forbidden subgraphs, Combinatorica $21$ $(2001)$ $155-170$.
\bibitem{kgg}K. Choromanski, Excluding pairs of tournaments, J. Graph Theory $89 $ $(2018)$ $ 266-287$.
\bibitem{polll}E. Berger, K. Choromanski, M. Chudnovsky, Forcing large transitive subtournaments, J. Comb. Theory, Ser. B. $112 $ $(2015)$ $ 1-17$.  
\bibitem{jhp}P. Erd\"{o}s and A. Hajnal, Ramsey-type theorems, Discrete Applied Mathematics $ 25$ $(1989)$ $ 37-52$.
\bibitem{bnmm}E. Berger, K. Choromanski, M. Chudnovsky, On the Erd\"{o}s$ - $Hajnal conjecture for six$ - $vertex tournaments, European Journal of Combinatorics $75$ $(2019)$ $113-122$.
\bibitem{kg}K. Choromanski, EH$-$suprema of tournaments with no nontrivial homogeneous sets, J. Comb. Theory, Ser. B. $114 $ $(2015)$ $ 97-123$.
\bibitem{ssss}S. Zayat, S. Ghazal, About the Erd\"{o}s-Hajnal conjecture for seven-vertex tournaments, submitted, arXiv:2010. 12331v1. 
\bibitem{sss}S. Zayat, S. Ghazal, Erd\"{o}s-Hajnal Conjecture for New Infinite Families of Tournaments, submitted, arXiv:2010. 12329v1.
\bibitem{ml} R. Stearns, The voting problem, Amer. Math. Monthly $66$ $ (1959)$ $ 761-763$.
\end{thebibliography}
\end{document}